\documentclass[margincite,leqno]{sl2art}

\usepackage{mathtools}

\usepackage{todonotes}

\usepackage{mathrsfs}

\usepackage{tikz,tikz-cd}
\usetikzlibrary{arrows.meta}
\tikzset{>={Stealth}}

\ExplSyntaxOn
\DeclareExpandableDocumentCommand{\IfEmptyTF}{mmm}
  {
    \tl_if_empty:nTF {#1}
      {#2}
      {#3}
  }
\ExplSyntaxOff

\def\leftfun#1{\mathopen{}\left#1}
\def\rightfun#1{\right#1}

\NewDocumentCommand{\ParenIfNonempty}{m}{%
  \IfEmptyTF{#1}{%
  }{%
    \leftfun(#1\rightfun)%
  }
}

\newcommand{\set}[1]{\left\{ \def\given{\ \middle| \ }  #1 \right\}  }
\newcommand{\defeq}{\mathrel{:=}}

\RenewDocumentCommand{\hom}{}{\operatorname{Hom}}

\NewDocumentCommand{\nr}{}{N}
\NewDocumentCommand{\pr}{}{p}
\NewDocumentCommand{\ii}{}{\mathsf{i}}
\NewDocumentCommand{\tu}{}{2 \pi \ii}

\NewDocumentCommand{\dif}{}{\,\mathrm{d}}

\NewDocumentCommand{\slg}{}{\operatorname{SL}_2(\mathbb{C})}

\NewDocumentCommand{\lN}{}{\operatorname{N}}
\NewDocumentCommand{\lS}{}{\operatorname{S}}
\NewDocumentCommand{\lE}{}{\operatorname{E}}
\NewDocumentCommand{\lW}{}{\operatorname{W}}

\NewDocumentCommand{\nrdel}{O{\nr}m}{\delta_{#1}\ParenIfNonempty{#2}}

\NewDocumentCommand{\pf}{O{} m}{\varphi_{#1}\IfEmptyTF{#2}{}{\leftfun(#2\rightfun)}}
\NewDocumentCommand{\psif}{O{} m}{\varphi_{#1}\IfEmptyTF{#2}{}{\leftfun(#2\rightfun)}}
\NewDocumentCommand{\psift}{O{} m}{\widetilde{\psi}_{#1}\IfEmptyTF{#2}{}{\leftfun(#2\rightfun)}}
\NewDocumentCommand{\pft}{O{} m}{\widetilde{\varphi}_{#1}\IfEmptyTF{#2}{}{\leftfun(#2\rightfun)}}
\NewDocumentCommand{\df}{O{} m }{\operatorname{D}_{#1}\ParenIfNonempty{#2}}
\NewDocumentCommand{\ldil}{m}{\operatorname{L}\ParenIfNonempty{#1}}

\NewDocumentCommand{\gauss}{m}{\Gamma\leftfun(#1\rightfun)}

\NewDocumentCommand{\wfunc}{m m}{w\leftfun(#1\middle|#2\rightfun)}

\NewDocumentCommand{\poch}{m O{q} m}{\left(#1; #2\right)_{#3}}

\NewDocumentCommand{\psisum}{m m m}{\Psi\leftfun( \begin{matrix} #1 \\ #2 \end{matrix}\middle| #3 \rightfun)}

\NewDocumentCommand{\U}{}{\mathsf{U}}
\NewDocumentCommand{\V}{}{\mathsf{V}}
\NewDocumentCommand{\Q}{}{\mathsf{Q}}
\NewDocumentCommand{\K}{}{\mathsf{K}}

\NewDocumentCommand{\qtor}{O{\omega}}{\mathcal{T}_{#1}}
\NewDocumentCommand{\qpmod}{m}{\mathcal{E}\ParenIfNonempty{#1}}
\NewDocumentCommand{\holbraid}{}{\mathsf{B}}
\NewDocumentCommand{\holker}{O{\nr} m m}{%
  \Omega_{#1}%
  \IfEmptyTF{#2}{%
    \IfEmptyTF{#3}{}{\leftfun(#3\rightfun)}%
  }{%
    \leftfun(#2%
    \IfEmptyTF{#3}{%
      \rightfun)%
    }{%
      \middle|#3\rightfun)%
    }%
  }
}
\NewDocumentCommand{\holkerhat}{O{\nr} m m}{%
  \widehat{\Omega}_{#1}%
  \IfEmptyTF{#2}{%
    \IfEmptyTF{#3}{}{\leftfun(#3\rightfun)}%
  }{%
    \leftfun(#2%
    \IfEmptyTF{#3}{%
      \rightfun)%
    }{%
      \middle|#3\rightfun)%
    }%
  }
}

\newcommand{\br}{to [out=00,in=180]}

\makeatletter
\def\qed@warning{}
\makeatother

\numberwithin{equation}{section}

\usepackage[export]{adjustbox}
\usepackage{graphicx,subcaption}
\graphicspath{{fig/}}

\usepackage{enumitem}

\newlist{thmenum}{enumerate}{1}
\setlist[thmenum]{label=(\alph*)}

\newkeytheorem{proposition}[style=theorem, parent=section]
\newkeytheorem{theorem, lemma, claim, corollary}[style=theorem,
  sibling=proposition]
\newkeytheorem{definition, remark, example}[style=definition,
  sibling=proposition]
\usepackage{zref-clever}

\NewDocumentCommand{\cref}{m}{\zcref{#1}}
\NewDocumentCommand{\Cref}{m}{\zcref[S]{#1}}

\renewbibmacro*{in:}{}
\title{
  A lift of Baxter's cyclic summation identity with applications to quantum topology
}
\author{Calvin McPhail-Snyder}
\address{Department of Mathematics, Duke University}
\email{calvin@sl2.site}


\begin{document}

\begin{abstract}
  While studying the chiral Potts model Baxter discovered a cyclic analogue of Ramanujan's \({}_1\psi_{1}\) summation identity.
  This identity was subsequently interpreted by Faddeev and Kashaev as a pentagon relation for the quantum dilogarithm function, a key ingredient in cyclic quantizations of Teichm\"uller theory.
  Existing forms of the identity involve multivalued functions so one must either accept overall phase ambiguities or keep track of branches on a complicated covering space.
  In this paper we show how to remove the phase ambiguities with a much simpler choice of covering space.
  As an application we compute the Fourier transform of the holonomy \(R\)-matrices of the author and Reshetikhin.
  A corollary is a simpler proof that Kashaev's \(R\)-matrix is a specialization of the colored Jones \(R\)-matrix, which was first shown by Murakami and Murakami.
\end{abstract}

\maketitle

\section{Introduction}

Consider the \defemph{shifted \(q\)-factorial}
\begin{equation}
  \poch{x}{n} = (1 - x) (1 - x q) \cdots (1 - x q^{n-1}).
\end{equation}
Here we focus on the \defemph{cyclic} case where \(q = \omega = e^{\tu / \nr}\) is an \(\nr\)th root of unity.
Because
\begin{equation*}
  \poch{x}[\omega]{\nr}
  =
  (1 - x)(1 - \omega x) \cdots (1 - \omega^{\nr-1} x)
  = 1 - x^{\nr}
\end{equation*}
we have quasi-periodic behavior
\begin{equation*}
  \poch{x}[\omega]{k + \ell \nr}
  =
  \poch{x}[\omega]{k} (1 - x^{\nr})^{\ell}.
\end{equation*}
One approach to handling this is to consider a function like
\begin{equation}
  \wfunc{x}{n} =
  \frac{(1 - x^{\nr})^{n/\nr}}{(1 - \omega x) (1 - \omega^{2} x) \cdots (1 - \omega^{n} x)}
\end{equation}
which is periodic: it depends only on \(n \mod \nr\).
The cost of imposing periodicity is keeping track of the branch of \((1 - x^{\nr})^{1/\nr}\).

The function \(w\) can be viewed as a cyclic \(q\)-analogue or quantization of the dilogarithm function \(\operatorname{Li}_2(z) = \sum_{n=1}^{\infty} z^{n}/n^{2}\).
From this perspective it is desirable to consider \(w\) as a function of a single argument that can be promoted to an operator.
One way to do this is to consider the multivalued function
\begin{equation*}
  \operatorname{d}_{\nr}(x) = \prod_{k=1}^{\nr-1} (1 - \omega^{k} x)^{(\nr-k)/\nr}
\end{equation*}
which satisfies (up to some power of \(\omega\)) the recurrence
\begin{equation*}
  \operatorname{d}_{\nr}(\omega x) = \frac{(1 - x^{\nr})^{1/\nr}}{1 - \omega x} \operatorname{d}_{\nr}(x)
  \text{ hence }
  \operatorname{d}_{\nr}(\omega^{k} x)
  =
  \wfunc{x}{k} \operatorname{d}_{\nr}(x).
\end{equation*}

In the study of the chiral Potts model \textcite{Baxter1988} discovered a cyclic analogue of Ramanujan's \({}_1\psi_{1}\) identity~\cite{Warnaar2012} for the shifted factorials \(\poch{x}[\omega]{n}\).
It was later interpreted as a \defemph{pentagon relation} for the \defemph{quantum dilogarithm} by \textcite{FaddeevKashaev1994}.
In terms of the function \(\operatorname{d}_{\nr}\) and parameters \(x, y, z\) satisfying
\begin{equation*}
  z^{\nr} = \frac{1 - x^{\nr}}{1 - y^{\nr}}
\end{equation*}
the identity is~\cite[equation (35)]{AuYang1996}
\begin{equation}
  \label{eq:cyclic-dilogarithm-summation}
  \sum_{k=0}^{\nr-1}
  \frac{\operatorname{d}_{\nr}(\omega^{k} x)}{\operatorname{d}_{\nr}(\omega^{k} y)}
  =
  \omega^{?} \frac{1 - z^{\nr}}{1 - z}
  \frac{\operatorname{d}_{\nr}(1)\operatorname{d}_{\nr}(yz/x)}
  {\operatorname{d}_{\nr}(y/x)\operatorname{d}_{\nr}(z)}
\end{equation}
for some power of \(\omega\).
Because of the \(\nr\)th roots in the definition of \(\operatorname{d}_{\nr}\) one either has to keep track of the branches or accept phase ambiguities.
In statistical mechanics one is frequently only concerned with the magnitude of the partition function so this is not a major problem.
However in quantum topology the partition function becomes a link invariant so these ambiguities are more significant: they lead to indeterminacies in quantum hyperbolic invariants of links and \(3\)-manifolds~\cite{BaseilhacBenedetti2003}.

\subsection{The main result}

We show how to resolve the phase ambiguities by using a different version of the quantum dilogarithm.
Consider a normalization of Faddeev's noncompact quantum dilogarithm~\cite{Faddeev2001} (see \cref{sec:dilogarithm-quantum}) satisfying%
\note{
  This \(\psif{}\) differs from the \(\varphi\) in~\cite{arXiv:2509.02365} by a constant.
}
\begin{equation*}
  \psif{\zeta + k} =
  \frac{\psif{\zeta}}{(1 - \omega^{\zeta+1}) \cdots (1 - \omega^{\zeta+k})}
  \quad
  \text{ and }
  \quad
  \psif{0} = 1
  .
\end{equation*}
For complex parameters \(\alpha, \beta, \gamma\) satisfying
\begin{equation}
  \label{eq:abg-relation}
  \frac{1 - \omega^{\nr \alpha}}{1 - \omega^{\nr \beta}} = \omega^{\nr \gamma}
\end{equation}
the function
\begin{equation*}
  k \mapsto \frac{\psif{\alpha + k}}{\psif{\beta + k}} \omega^{k \gamma}
\end{equation*}
has period \(\nr\) and the sum
\begin{equation}
  \psisum{\alpha}{\beta}{\gamma} \defeq \sum_{k=0}^{\nr-1} \frac{\psif{\alpha + k}}{\psif{\beta + k}} \omega^{k \gamma}
\end{equation}
is a cyclic version of a \({}_1\psi_{1}\) sum.

\begin{theorem}
  \label{thm:1psi1-summation}
  For parameters satisfying~\eqref{eq:abg-relation}
  \begin{align}
    \label{eq:1psi1-summation} \psisum{\alpha}{\beta}{\gamma}         & = \nr \omega^{(\nr-1-\beta)\gamma} \frac{\psif{-\gamma}\psif{\alpha - \beta + \nr - 1}}
    {\psif{\alpha - \beta - \gamma + \nr - 1}}, \\
    \label{eq:1psi1-summation-flipped} \psisum{\alpha}{\beta}{\gamma} & = \omega^{-\alpha\gamma} \frac{\psif{\gamma + \beta - \alpha}}
    {\psif{\gamma + \nr - 1}\psif{\beta - \alpha}}.
  \end{align}
  In particular both sides are meromorphic functions.
\end{theorem}
While the identities above require a choice of \(\alpha, \beta, \gamma\) they depend on these choices in a simple way.
The underlying parameters \(a = e^{\tu \alpha}, b = e^{\tu \beta}, c = e^{\tu \gamma}\) lie in the set
\begin{equation*}
  \mathcal{Y} =
  \set{ a, b, c \ne 0 \given c = \frac{1 -a}{1-b} }
\end{equation*}
and the identities in \cref{thm:1psi1-summation} hold on a covering space of \(\mathcal{Y}\) with abelian deck transformation group \(\mathbb{Z}^{3}\).
In contrast the covering space underlying the summation identity~\eqref{eq:cyclic-dilogarithm-summation} is much more complicated.
For the identities of \cref{thm:1psi1-summation} one wants to view \(\operatorname{d}_{\nr}(x)\) as a function of \(x^{\nr}\), in which case it is defined on a covering space of \(\mathbb{C} \setminus \set{0,1}\) with a nonabelian group of deck transformations.
When considering multiple parameters like in~\eqref{eq:cyclic-dilogarithm-summation} these combine in a complicated way.

Some simpler, closely related identities are the quasi-periodic discrete Fourier transforms:
\begin{theorem}
  \label{thm:single-qlog-fourier}
  If \(\beta, \gamma\) satisfy \(\omega^{\nr \gamma} = (1 - \omega^{\nr \beta})^{-1}\) then
  \begin{align*}
    \sum_{k=0}^{\nr-1} \frac{\omega^{(\beta+k)\gamma}}{\psif{\beta + k}} & = \frac{\nr}{c_{\nr}} \omega^{(\nr-1)\gamma}\psif{-\gamma}, \\
    \sum_{k=0}^{\nr-1} \psif{\beta + k}\omega^{-\gamma(k+\beta)}         & = \frac{c_{\nr}}{\psif{-\gamma + \nr - 1}}
  \end{align*}
  where
  \[
    c_{\nr}
    =
    \nr^{-1/2}\exp\frac{\pi\ii}{12}
    \left( \nr - 3 + \frac{1}{\nr} \right).
    \qedhere
  \]
\end{theorem}

Our proof follows the algebraic approach to the five-term relation introduced by \textcite{FaddeevKashaev1994}.
Simple modules \(\qpmod{a,b}\) over the quantum torus \(\qtor = \langle \U, \V | \U \V = \omega \V \U \rangle\) are classified by their central characters \(\chi(\U^{\nr}) = a, \chi(\V^{\nr}) = b\).
We define a mutation operator \(\K \colon \qpmod{a,b} \to \qpmod{b(1 - (-1)^{\nr-1} a), a^{-1}}\) via \(\psif{}\) and a discrete Fourier transform.
\(\K^{5}\) is a \(\qtor\)-module endomorphism of \(\qpmod{a,b}\) so it acts by a scalar.
We determine that scalar using determinants and a classical five-term identity; the reduction to the classical dilogarithm comes from an identity of \textcite{Garoufalidis2014}.
The resulting matrix identity produces a more complicated identity~\eqref{eq:five-psi-summation} that gives \cref{thm:single-qlog-fourier} as a limiting case, and we can then extract \cref{eq:1psi1-summation}.

In some sense this is the ``same'' as Baxter's original proof of \cref{eq:cyclic-dilogarithm-summation}.
His technique (as described by \textcite{AuYang1996}) is to compute the determinant of a circulant matrix in two ways.
By factoring one side and counting degrees one can establish that the identity is true up to a constant function on parameter space, which can then be found by a limit.
A similar argument (not phrased in terms of determinants) is given by \textcite[Appendix C]{Kashaev1993}.
In our proof we follow \textcite{FaddeevKashaev1994} and use representation theory to avoid the somewhat complicated computations involved in factoring the determinant.

\subsection{Application to holonomy \texorpdfstring{\(R\)}{R}-matrices}

In recent work~\cite{arXiv:2509.02354} the author and Reshetikhin found a family of \defemph{holonomy \(R\)-matrices} pa\-ra\-me\-trized by the space of flat \(\mathfrak{sl}_2\) connections on the complement of a braid generator; there is also an interpretation in terms of the geometry of ideal octahedra~\cite{McPhailSnyder2022,McphailSnyder2024,arXiv:2509.02365}.
These matrices can be viewed as a geometric deformation of Kashaev's \(R\)-matrix~\cite{Kashaev1995}, and the corresponding link invariants can be interpreted as a quantization of the \(\slg\) Chern-Simons invariant~\cite{arXiv:2509.02365}.
The \(R\)-matrices act on the modules \(\qpmod{a,b}\) and are defined by their matrix coefficients in a certain basis.
In \cref{sec:R-matrix} we compute the matrix coefficients in the Fourier dual basis.
\Cref{thm:holonomy-fourier} generalizes our earlier result~\cite[Theorem 5.9]{arXiv:2509.02354}, which only computed these coefficients in a geometrically degenerate limit.
By using \cref{thm:1psi1-summation} we get a more general result with a simpler proof.
A corollary is a simpler proof of a result of \textcite[Proposition 4.2]{Murakami2001}: Kashaev's \(R\)-matrix is the discrete Fourier transform of the colored Jones \(R\)-matrix at a root of unity up to some elementary factors.

\subsection{Related work and future directions}

\textcite{IpYamazaki2016} give a general method for obtaining cyclic dilogarithm identities from cycles of mutations.
In that language our identity comes from the \(A_{2}\) quiver, as in~\cite[Example 3.12]{IpYamazaki2016}.
The new contribution of this paper is to turn these into identities of meromorphic functions.
This method should also work for identities obtained from other cluster mutation sequences.

Recent work of \textcite{Ishibashi2025} systematically constructs \defemph{cyclic quantum Teichm\"uller theory} using cyclic dilogarithms related to \(\operatorname{d}_{\nr}\).
The resulting mapping class group representations have phase ambiguities unless one keeps careful track of branches~\cite{Ishibashi2026}.
Our lifted identity may provide a different way to handle these phases: by replacing the modules associated to edges with the spaces \(\qpmod{a,b}\) of quasi-periodic functions and the dilogarithms with \(\psif{}\) one should be able to remove the phase ambiguity in the pentagon relation.

\subsection*{Notation}

We write \(\ii\) for the imaginary unit.
Throughout \(\nr \ge 2\) is an integer and \(\omega^{x} = e^{\tu x / \nr}\).
In \cref{sec:quantum,sec:R-matrix} we will abbreviate \(\pr = (\nr - 1)/2\) and \(\epsilon = \omega^{\nr p} = (-1)^{\nr-1}\).

\subsection*{Acknowledgements}

I would like to thank Rei Inoue for helpful conversations.

\subsection*{AI use statement}

I developed proofs of the new results (\cref{thm:1psi1-summation,thm:holonomy-fourier}) independently of AI assistance.
However these proofs were more complicated than the ones appearing in the text and were not self-contained: the proof of \cref{thm:1psi1-summation} required Baxter's identity~\eqref{eq:cyclic-dilogarithm-summation} as an input.
Via AI assistance I found the simpler, better proofs presented here.
Most of the contribution of the AI agent was to check ideas and search the literature much faster than I could unassisted, although one new idea (to me) was the definition of the mutation \(\K\) using a Fourier transform instead of a derivative operator.
This avoids some technical complications (see \cref{sec:mutations} for details).
Another AI contribution was taking boundary limits on parameter space to evaluate constants, like in the last part of the proof of \cref{thm:weyl-pentagon-mutation}.
These computations are necessary but are typically fiddly and not very interesting, and for that reason are frequently omitted.
I have decided to include these details since I did not have to write them out manually; the proofs in the text were found via AI and then cleaned up.

\section{Facts about dilogarithms}
\label{sec:dilogarithms}

All of the results here are well known.
However some of the proofs are folklore, and in particular I am not aware of a published source for the residue computations in the proofs of \cref{thm:psi-recurrence,thm:psi-inversion}, although I am confident they were known to \textcite{Faddeev2001}.
Similarly \cref{thm:exact-quantum-dilogarithm} and its proof are a special case of a theorem of \textcite[Theorem 1.9]{Garoufalidis2014}.

\subsection{The classical dilogarithm}
\label{sec:dilogarithm-classical}

The power series
\begin{align*}
  - \log(1 - z) & = \sum_{k=1}^{\infty} \frac{z^{k}}{k} & \operatorname{Li}_{2}(z) & = \sum_{k=1}^{\infty} \frac{z^{k}}{k^{2}}
\end{align*}
extend to analytic functions on \(\mathbb{C} \setminus [1, \infty)\), so for \(\Im \zeta > 0\) we may define
\begin{equation}
  \label{eq:ldil-def}
  \ldil{\zeta}
  = \frac{\operatorname{Li}_2(e^{\tu\zeta})}{\tu}
  + \zeta\log(1 - e^{\tu\zeta})
  + \frac{\pi \ii}{12}
  .
\end{equation}

\begin{proposition}
  \label{thm:classical-dilogarithm-continuation}
  The function \(e^{\ldil{\zeta}}\) extends meromorphically to \(\mathbb{C}\) and obeys the relation
  \begin{equation}
    \label{eq:classical-dilogarithm-recurrence}
    e^{\ldil{\zeta + 1}}
    =
    (1 - e^{\tu\zeta})
    e^{\ldil{\zeta}}
    .
  \end{equation}
  It has zeros at positive integers and poles at negative integers.
  In particular it is holomorphic near \(0\) with \(\ldil{0} = 0\).
  It satisfies the \defemph{reflection identity}
  \begin{equation}
    \label{eq:inversion-classical}
    \exp \leftfun( \ldil{\zeta} + \ldil{-\zeta} \rightfun)
    =
    e^{\pi\ii\zeta^2}.
  \end{equation}
\end{proposition}

\begin{proof}
  Write \(\ell(\zeta)\) for the derivative of \(\ldil{\zeta}\).
  On the upper half-plane we have
  \[
    \ell(\zeta)
    =
    \zeta
    \frac{\dif}{\dif \zeta}
    \log(1 - e^{\tu \zeta})
    =
    -\tu\zeta\frac{e^{\tu\zeta}}{1 - e^{\tu\zeta}}
    .
  \]
  Thus \(\ell\) is meromorphic on \(\mathbb{C}\), and near each integer \(m\) we have
  \[
    \ell(\zeta) = \frac{m}{\zeta - m} + h_m(\zeta)
  \]
  where \(h_m\) is holomorphic.
  In particular the singularity at \(0\) is removable.

  We can analytically continue \(\ldil{}\) in \(\mathbb{C} \setminus \mathbb{Z}\) by integrating \(\ell\).
  The residues of \(\ell\) are all integers, so \(e^{\ldil{}}\) is single-valued.
  To check the value at \(0\) observe that as \(\zeta \to 0\) from the upper half-plane \(\zeta\log(1 - e^{\tu\zeta}) \to 0\) and \(\operatorname{Li}_2(e^{\tu\zeta}) \to \pi^2/6\).

  Similar reasoning gives the reflection relation: we have
  \begin{equation*}
    \frac{\mathrm d}{\mathrm d\zeta}
    \left( \ldil{\zeta} + \ldil{-\zeta} \right)
    =
    \ell(\zeta) - \ell(-\zeta)
    = \tu\zeta.
  \end{equation*}
  Both terms tend to \(0\) as \(\zeta \to 0\), so we conclude that
  \begin{equation*}
    \ldil{\zeta} + \ldil{-\zeta} = \pi \ii \zeta^{2}
  \end{equation*}
  in a neighborhood of \(0\), and the identity follows by analytic continuation.
\end{proof}

\textcite[Proposition 2.5]{Neumann2004} established a lifted five-term relation for the Rogers dilogarithm extending the classical one for \(\operatorname{Li}_2\)~\cite{Zagier2007}.
The identity below for our version \(\ldil{}\) is equivalent.
The proof is well known, but we give it in detail to emphasize the similarities with the proof of \cref{thm:weyl-pentagon-mutation}.

\begin{theorem}
  \label{thm:five-term-classical}
  Consider the set \(\mathcal{Z} \subset \mathbb{C}^{\mathbb{Z}/5\mathbb{Z}}\) of points \(z_{i}\) indexed cyclically mod \(5\) satisfying
  \begin{equation}
    \label{eq:five-term-classical-shapes}
    1 - z_i
    = z_{i-1}z_{i+1}.
  \end{equation}
  For any logarithms \(\zeta_{i}\) satisfying \(e^{\tu \zeta_{i}} = z_{i}\) we have the
  \defemph{symmetric five-term relation}
  \begin{equation}
    \label{eq:five-term-classical-cyclic}
    \exp\leftfun(
    \sum_{i=0}^{4}\ldil{\zeta_i}
    - \tu\sum_{i=0}^{4}\zeta_i\zeta_{i+1}
    \rightfun)
    = e^{\pi\ii/6}.
  \end{equation}
\end{theorem}

\begin{proof}
  We first show that the identity is really a fact about functions on \(\mathcal{Z}\): it does not depend on the choice of logarithms \(\zeta_{i}\).
  Writing
  \begin{equation*}
    \Lambda(\zeta_{0}, \dots, \zeta_{4})
    =
    \sum_{i=0}^{4}\ldil{\zeta_i}
    - \tu\sum_{i=0}^{4}\zeta_i\zeta_{i+1}
  \end{equation*}
  \cref{eq:classical-dilogarithm-recurrence} gives
  \begin{align*}
    \Lambda(\dots, \zeta_{i} + n_i, \dots)
    -
    \Lambda(\zeta_{0}, \dots, \zeta_{4})
    =
    n_{i}\log(1 - z_{i}) + \tu k_{i}
    - \tu n_{i}(\zeta_{i-1} + \zeta_{i+1})
  \end{align*}
  where \(k_{i}\) is an integer accounting for any discrepancies in the choice of branch.
  By \cref{eq:five-term-classical-shapes},
  \begin{equation*}
    \log(1 - z_{i}) \equiv \tu ( \zeta_{i-1} + \zeta_{i+1} ) \pmod{\tu \mathbb{Z}}.
  \end{equation*}
  Thus \(\exp \leftfun( \Lambda(\dots, \zeta_{i} + n_i, \dots) - \Lambda(\zeta_{0}, \dots, \zeta_{4}) \rightfun) = 1\).
  By successively changing by integers we can move between any two choices of \((\zeta_{i})_{i}\) so the value of the left-hand side of \cref{eq:five-term-classical-cyclic} is independent of this choice as claimed.

  Furthermore the set \(\mathcal{Z}\) is connected: \(z_0,z_1\) determine the remaining shapes by
  \begin{equation*}
    z_2 = \frac{1 - z_1}{z_0},
    \qquad
    z_3 = \frac{z_0 + z_1 - 1}{z_0z_1},
    \qquad
    z_4 = \frac{1 - z_0}{z_1},
  \end{equation*}
  and the allowed \((z_0,z_1)\) form the complement of finitely many lines in \(\mathbb C^2\).
  Using analytic continuation we can reduce the proof to checking the identity on the \defemph{positive real chamber}
  \begin{equation*}
    \mathfrak{C} =
    \set{x \in \mathcal{Z} \given 
0 < x_{0}, x_{1} < 1
\text{ and }
x_{0} + x_{1} > 1
}.
  \end{equation*}
  On \(\mathfrak{C}\) we have \(x_i \in (0,1)\) so with the principal branch of \(\log\) and \(\xi_{i} = \frac{1}{\tu}\log x_{i}\),
  \begin{equation*}
    \log(1 - e^{\tu\xi_i})
    = \tu(\xi_{i-1} + \xi_{i+1})
    .
  \end{equation*}
  Since \(\mathrm d\ldil{\zeta} = \zeta\,\mathrm d\log(1 - e^{\tu\zeta})\) we see that
  \begin{align*}
    \mathrm d\sum_{i=0}^{4}\ldil{\xi_i} = \tu\sum_{i=0}^{4}\xi_i (\mathrm d\xi_{i-1} + \mathrm d\xi_{i+1})
    = \tu\,\mathrm d\sum_{i=0}^{4}\xi_i\xi_{i+1}.
  \end{align*}
  Thus \(\mathrm{d} \Lambda = 0\), so \(\Lambda\) is a locally constant function on \(\mathfrak{C}\), and because \(\mathfrak{C}\) is connected it is constant.

  To evaluate the constant set \(x_0 = x_1 = t\) for \(1/2 < t < 1\).
  The other parameters are then
  \begin{equation*}
    x_2 = x_4 = \frac{1 - t}{t},
    \qquad
    x_3 = \frac{2t - 1}{t^2}.
  \end{equation*}
  As \(t \to 1\) from below \(x_0,x_1,x_3 \to 1\) and \(x_2,x_4 \to 0\), so every product \(\xi_i\xi_{i+1}\) tends to zero since the two parameters tending to \(\ii\infty\) are not adjacent.
  We conclude that
  \begin{equation*}
    \Lambda(\xi_{0}, \dots, \xi_{4})
    =
    \sum_{i=0}^{4}\ldil{\xi_i}
    - \tu\sum_{i=0}^{4}\xi_i\xi_{i+1}
    \to 3 \ldil{0} + 2\lim_{t \to +\infty}\ldil{\ii t}
    = \frac{\pi \ii}{6}.
    \qedhere
  \end{equation*}
\end{proof}

\begin{corollary}
  \label{thm:classical-z0z1-inversion}
  Suppose \(\beta, \gamma\) satisfy \(e^{\tu \beta} + e^{\tu \gamma} = 1.\) Then
  \begin{equation}
    \label{eq:classical-complementary-inversion}
    \exp\leftfun( \ldil{\beta} + \ldil{\gamma} \rightfun)
    =
    \exp\leftfun( \frac{\pi\ii}{12} + \tu\beta\gamma \rightfun).
    \qedhere
  \end{equation}
\end{corollary}

\begin{proof}
  We first prove this in the positive real chamber \(\mathfrak{C}\) from the proof above.
  If we set \(z_1 = b \in (0,1)\) and let \(z_0 \to 1^-\) then the five shapes tend to \((1,b,1 - b,1,0)\).
  Writing \(\beta = \log b/\tu\) and \(\gamma = \log(1 - b)/\tu\) the quadratic sum tends to \(\beta\gamma\) and the dilogarithm terms tend to
  \begin{equation*}
    2\ldil{0} + \ldil{\beta} + \ldil{\gamma}
    + \lim_{t \to +\infty}\ldil{\ii t}
    =
    \ldil{\beta} + \ldil{\gamma} + \frac{\pi\ii}{12}.
  \end{equation*}
  Thus the limit of~\eqref{eq:five-term-classical-cyclic} gives~\eqref{eq:classical-complementary-inversion} for these logarithms.
  The general case now follows from analytic continuation; both sides of the relation transform in the same way when we change \(\beta\) or \(\gamma\).
\end{proof}

We can use \cref{thm:classical-z0z1-inversion} to rewrite the five-term relation~\eqref{eq:five-term-classical-cyclic} as
\begin{equation}
  \label{eq:five-term-classical}
  \frac{\exp \leftfun( \ldil{\beta} + \ldil{-\gamma} + \ldil{\alpha - \beta} \rightfun)}{\exp \leftfun( \ldil{\alpha} + \ldil{\alpha - \beta -\gamma} \rightfun)}
  =
  e^{-\tu\beta\gamma}
\end{equation}
which corresponds directly to the Baxter identity~\eqref{eq:1psi1-summation-flipped}.
Here the parameters \(\alpha, \beta, \gamma\) satisfy~\eqref{eq:abg-relation}, that is,
\begin{equation*}
  c = e^{\tu \gamma} = \frac{1 - a}{1 - b}
\end{equation*}
for \(a = e^{\tu \alpha}\), \(b = e^{\tu \beta}\), and \(c = e^{\tu \gamma}\).

To obtain~\eqref{eq:five-term-classical} consider the shapes
\begin{equation*}
  (z_0,z_1,z_2,z_3,z_4)
  = \left( \frac{a}{b},1 - a,b, c^{-1}, 1 - \frac{a}{bc} \right).
\end{equation*}
Choose logarithms
\begin{equation*}
  \zeta_0 = \alpha - \beta,
  \qquad
  \zeta_2 = \beta,
  \qquad
  \zeta_3 = -\gamma,
\end{equation*}
which have \(\zeta_0 + \zeta_2 = \alpha\) and \(\zeta_0 + \zeta_3 = \alpha - \beta - \gamma\).
Now \cref{thm:five-term-classical} says that
\begin{equation*}
  \ldil{\alpha - \beta}
  +
  \ldil{\zeta_{1}}
  +
  \ldil{\beta}
  +
  \ldil{-\gamma}
  +
  \ldil{\zeta_{4}}
  +
  \zeta_1(\zeta_0 + \zeta_2)
  + \zeta_4(\zeta_0 + \zeta_3) - \tu \beta\gamma
\end{equation*}
is equal to \(\pi \ii / 6\) modulo \(\tu \mathbb{Z}\).
Simplifying the left-hand side using~\eqref{eq:classical-complementary-inversion} on the pairs \((\zeta_1,\zeta_0 + \zeta_2)\) and \((\zeta_4,\zeta_0 + \zeta_3)\) gives the identity.

\subsection{The quantum dilogarithm}
\label{sec:dilogarithm-quantum}

Our function \(\psif{}\) is a version of Faddeev's noncompact quantum dilogarithm~\cite[Section 6,
  equation (15)]{Faddeev2001}.
Writing their \(e_{b}\) in the more common notation \(\Phi_{\mathsf b}\) we have
\begin{equation}
  \label{eq:psi-faddeev}
  \psif[\nr]{\zeta}
  =
  \Phi_{\sqrt\nr}\left( \frac{\ii\zeta}{\sqrt\nr} + \frac{\ii(1 - \nr)}{2\sqrt\nr} \right)
  \Phi_{\sqrt\nr}\left( \frac{\ii(1 - \nr)}{2\sqrt\nr} \right)^{-1}
\end{equation}
so that \(\psif{0} = 1\).
We give a self-contained definition:

\begin{definition}
  \label{thm:normalized-quantum-dilogarithm}
  For \(-1 < \Re\zeta < \nr\), define
  \begin{equation}
    \label{eq:psi-integral}
    \psif[\nr]{\zeta}
    \defeq
    \exp\left[ \int_{\mathbb{R}+\ii\epsilon} \frac{e^{(2\zeta-\nr+1)w/\nr} - e^{(1-\nr)w/\nr}}
      {4w\sinh w\sinh(w/\nr)}\,\mathrm dw \right],
  \end{equation}
  where \(+ \ii \epsilon\) means we detour above the pole at \(0\) as in \cref{fig:contour-completion}.
  We usually abbreviate \(\psif{\zeta} = \psif[\nr]{\zeta}\).
\end{definition}

\begin{proposition}[\cite{Faddeev2001}]
  \label{thm:psi-recurrence}
  The function \(\psif{}\) extends meromorphically to \(\mathbb{C}\) via the functional equation
  \begin{equation}
    \label{eq:psi-recurrence}
    \psif{\zeta + 1}
    =
    \frac{\psif{\zeta}}{1 - \omega^{\zeta+1}}.
  \end{equation}
\end{proposition}

\begin{proof}
  Consider the overlap \(-1 < \Re\zeta < \nr - 1\) of the convergence regions for \(\psif{\zeta}\) and \(\psif{\zeta + 1}\).
  Subtracting the exponents in~\eqref{eq:psi-integral} gives
  \[
    \frac{\psif{\zeta + 1}}{\psif{\zeta}}
    =
    \exp\left[ \int_{\mathbb{R}+\ii\epsilon} \frac{e^{(2\zeta-\nr+2)w/\nr}}{2w\sinh w}\,\mathrm dw \right].
  \]
  Suppose also that \(\Im\zeta > 0\).
  By completing the contour as in \cref{fig:contour-completion} we can show that the integral is given by the sum of residues in the upper half-plane.
  \begin{figure}[tbp]
    \centering
    \begin{tikzpicture}[
    x=1cm,
    y=1cm,
    font=\small,
    original/.style={draw=accent, thick},
    completion/.style={draw=slblue, thick},
    pole/.style={draw=slred, thick}
  ]
  \draw[gray!65, ->] (-4.6, 0) -- (4.7, 0)
  node[right, text=black] {\(\Re w\)};
  \draw[gray!65, ->] (0, -0.55) -- (0, 3.75)
  node[above, text=black] {\(\Im w\)};

  \draw[original] (-4.5, 0) -- (-0.35, 0)
  arc[start angle=180, end angle=0, radius=0.35] -- (4.5, 0);
  \draw[original, ->] (-4.2, 0) -- (-3.8, 0);
  \draw[original, ->] (1.5, 0) -- (2.1, 0);
  \draw[original, ->] (3.8, 0) -- (4.2, 0);
  \node[text=accent, above] at (3.9, 0.1)
  {\(\mathbb{R} + \ii\epsilon\)};

  \draw[completion] (3, 0) -- (3, 3.15) -- (-3, 3.15) -- (-3, 0);
  \draw[completion, ->] (3, 1.2) -- (3, 1.85);
  \draw[completion, ->] (1.6, 3.15) -- (0.9, 3.15);
  \draw[completion, ->] (-3, 1.85) -- (-3, 1.2);
  \node[below] at (-3, -0.08) {\(-(M + 1/2)\pi\)};
  \node[below] at (3, -0.08) {\((M + 1/2)\pi\)};
  \node[above left, text=slblue] at (-3, 3.15)
  {\(\ii \Gamma_M\)};

  \foreach \h/\lab in {0/0, 0.9/{\pi\ii}, 1.8/{2\pi\ii}, 2.7/{M\pi\ii}}
  {
    \draw[pole] (-0.065, \h - 0.065) -- (0.065, \h + 0.065);
    \draw[pole] (-0.065, \h + 0.065) -- (0.065, \h - 0.065);
    \ifdim\h pt=0pt
    \node[below right, text=slred] at (0.12, -0.08) {\(\lab\)};
    \else
    \node[right, text=slred] at (0.12, \h) {\(\lab\)};
    \fi
  }
  \node[fill=white, inner sep=1pt] at (0, 2.25) {\(\vdots\)};
\end{tikzpicture}
    \caption{
      Integration contours used for \(\psif{}\).
      The original contour in \textcolor{accent}{gold} is along the real axis with a detour around \(0\).
      In the proofs of \cref{thm:psi-recurrence,thm:exact-quantum-dilogarithm} we complete it by adding the \textcolor{slblue}{blue} contour \(\Gamma_{M}\) and taking \(M \to \infty\).
    }
    \label{fig:contour-completion}
  \end{figure}
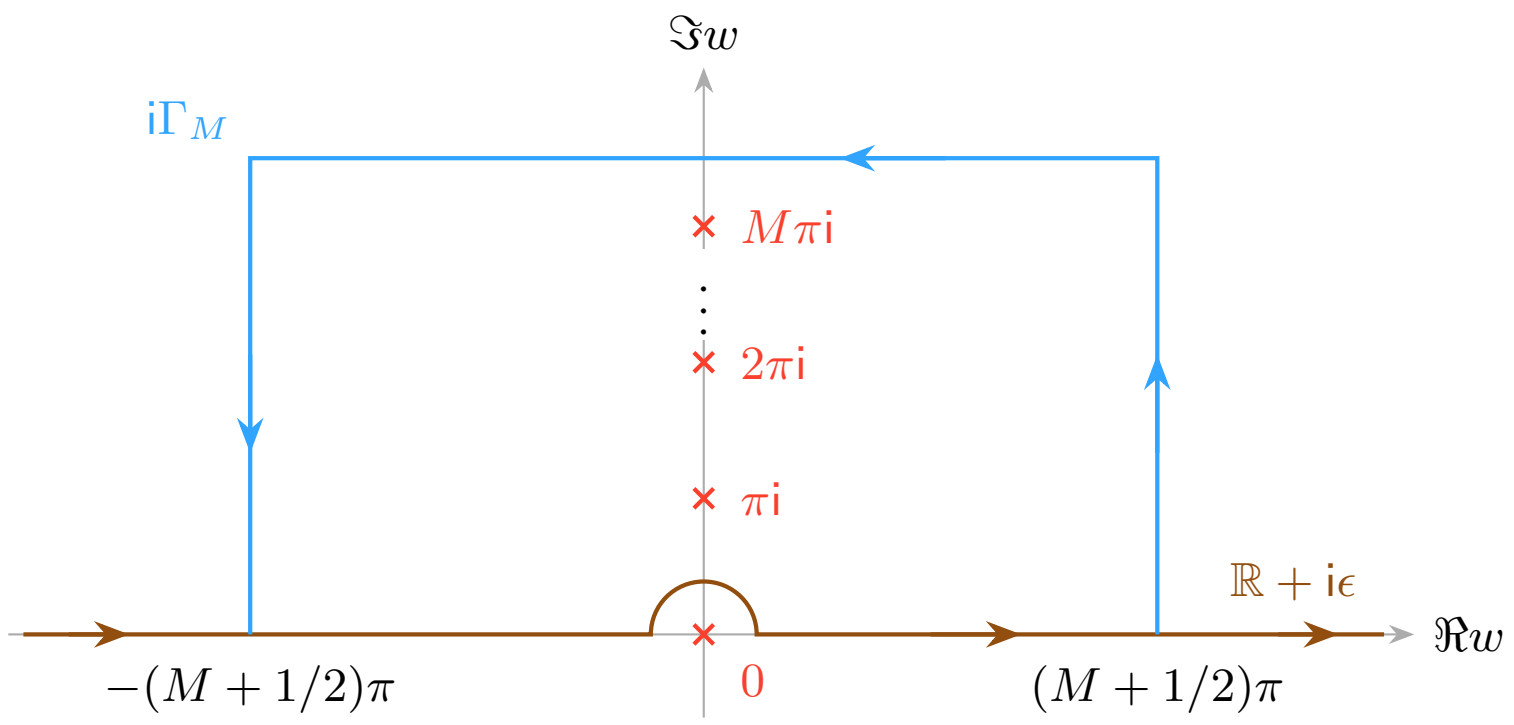
  This gives
  \[
    \int_{\mathbb{R}+\ii\epsilon}
    \frac{e^{(2\zeta-\nr+2)w/\nr}}{2w\sinh w}\,\mathrm dw
    =
    \sum_{k=1}^{\infty}\frac{\omega^{k(\zeta+1)}}{k}
    =
    -\log(1 - \omega^{\zeta+1})
  \]
  which gives the identity on our set; analytic continuation completes the proof.
\end{proof}

In the introduction we discussed another version of the quantum dilogarithm denoted \(\operatorname{d}_{\nr}\), which is a function of \(\omega^{\zeta}\).
To compare it with \(\psif{}\) we use the logarithmic version
\begin{equation*}
  \operatorname{D}_{\nr}(\zeta)
  =
  \exp \frac{1}{\nr}
  \sum_{j=1}^{\nr-1}(\nr - j)\log(1 - \omega^{\zeta+j})
\end{equation*}
which is multivalued because of the branch cuts of the logarithms.
For suitable branches they satisfy \(\operatorname{D}_{\nr}(\zeta) = \operatorname{d}_{\nr}(\omega^{\zeta})\) and thus
\begin{equation*}
  \operatorname{D}_{\nr}(\zeta + 1)
  =
  \frac{(1 - \omega^{\nr \zeta})^{1/\nr}}{1 - \omega^{\zeta+1}}
  \operatorname{D}_{\nr}(\zeta)
  .
\end{equation*}
\textcite[Theorem 1.9]{Garoufalidis2014} showed that one can remove the phase ambiguities by taking a suitable ratio of \(\operatorname{D}_{\nr}\) and \(e^{\ldil{}}\):
\begin{theorem}
  \label{thm:exact-quantum-dilogarithm}
  \begin{equation}
    \label{eq:exact-psi-value}
    \psif[\nr]{\zeta}
    =
    c_{\nr} \omega^{1/24}
    \operatorname{D}_{\nr}(\zeta)
    e^{-\ldil{\zeta}/\nr}
  \end{equation}
  It is not obvious that the right-hand side of \cref{eq:exact-psi-value} is a single-valued meromorphic function on \(\mathbb{C}\); this is part of the claim.
  The constant \(c_{\nr}\) is the same as in \cref{thm:single-qlog-fourier}:
  \begin{equation*}
    c_{\nr}
    =
    \nr^{-1/2}
    \exp\leftfun[ \frac{\pi\ii}{12}
    \left( \nr - 3 + \frac{1}{\nr} \right) \rightfun]
    \qedhere
  \end{equation*}
\end{theorem}

\begin{proof}
  This is a special case of~\cite[Theorem 1.9]{Garoufalidis2014}, but to make the paper more self-contained we give a proof.
  As before assume \(-1 < \Re\zeta < \nr\) and \(\Im\zeta > 0\).
  Set
  \[
    f_\zeta(w)
    =
    \frac{e^{(2\zeta-\nr+1)w/\nr}}{4w\sinh w\sinh(w/\nr)},
    \qquad
    I(\zeta)
    =
    \int_{\mathbb{R}+\ii\epsilon} f_\zeta(w)\,\mathrm dw
  \]
  so that \(\psif{\zeta} = \exp[I(\zeta) - I(0)]\).
  A contour argument similar to the one in the proof of \cref{thm:psi-recurrence} shows that
  \[
    I(\zeta)
    =
    \tu\sum_{k=1}^{\infty}\operatorname{Res}(f_\zeta,\pi\ii k).
  \]
  The simple poles and double poles have residues
  \begin{align*}
    \tu\operatorname{Res}(f_\zeta,\pi\ii k)    & = \frac{\omega^{k\zeta}}{k(1 - \omega^{-k})}                                            &  & (\nr \nmid k), \\
    \tu\operatorname{Res}(f_\zeta,\pi\ii\nr m) & = e^{\tu m\zeta}\left[ \frac{\zeta - (\nr - 1)/2}{\nr m} - \frac{1}{\tu\nr m^2} \right] &  & (m \ge 1).
  \end{align*}
  Both residue series converge absolutely.
  Differentiating the finite geometric sum gives
  \[
    \sum_{j=1}^{\nr-1}(\nr - j)\omega^{jk}
    =
    \begin{cases}
      -\nr/(1 - \omega^{-k}), & \nr \nmid k, \\
      \nr(\nr - 1)/2,         & \nr \mid k.
    \end{cases}
  \]
  Expanding the logarithms and dilogarithm into their power series and doing some algebra yields
  \begin{equation}
    \label{eq:unnormalized-dilogarithm-integral}
    I(\zeta)
    =
    \frac{1}{\nr}\left[ -\ldil{\zeta} + \sum_{j=1}^{\nr-1}(\nr - j)\log(1 - \omega^{\zeta+j}) + \frac{\pi\ii}{12} \right].
  \end{equation}
  The last term cancels the constant in~\eqref{eq:ldil-def}.

  To evaluate \(I(0)\) take \(\zeta \to 0\) from above.
  For \(1 \le j < \nr\), the principal branches give
  \[
    \log(1 - \omega^j)
    =
    \log\bigl(2\sin(\pi j/\nr)\bigr)
    + \pi\ii(j/\nr - 1/2).
  \]
  Pairing \(j\) with \(\nr - j\), using \(\prod_{j=1}^{\nr-1}(1 - \omega^j) = \nr\), and summing the imaginary parts gives
  \[
    \sum_{j=1}^{\nr-1}(\nr - j)\log(1 - \omega^j)
    =
    \frac{\nr}{2}\log\nr
    - \frac{\pi\ii}{12}(\nr - 1)(\nr - 2).
  \]
  Since \(\ldil{0} = 0\), we obtain
  \[
    I(0)
    =
    \frac{1}{2}\log\nr
    + \frac{\pi\ii}{12}\left( 3 - \nr - \frac{1}{\nr} \right).
  \]
  Thus \(e^{-I(0)} = c_{\nr}\), and exponentiating~\eqref{eq:unnormalized-dilogarithm-integral} gives
  \[
    \psif{\zeta}
    =
    e^{I(\zeta)-I(0)}
    =
    c_{\nr} \omega^{1/24}
    \operatorname{D}_{\nr}(\zeta)
    e^{-\ldil{\zeta}/\nr}.
  \]
  This proves~\eqref{eq:exact-psi-value} for \(-1 < \Re\zeta < \nr\) and \(\Im\zeta > 0\) and the general case follows by analytic continuation.
\end{proof}

We can use similar methods to find a quantum analogue of the reflection identity~\eqref{eq:inversion-classical}.
It is convenient to state it in terms of the \defemph{shifted Gaussian kernel}
\begin{equation*}
  \gauss{\zeta} = \exp\frac{\pi\ii}{\nr}\zeta(\zeta + \nr - 1)
\end{equation*}
which in particular satisfies
\begin{equation}
  \label{eq:gaussian-key-property}
  \gauss{\alpha+\beta} = \gauss{\alpha}\gauss{\beta}\omega^{\alpha\beta}
\end{equation}

\begin{proposition}
  \label{thm:psi-inversion}
  \begin{equation}
    \label{eq:psi-inversion}
    \psif{\zeta}\psif{\nr - 1 - \zeta}
    =
    \frac{1}{\nr\gauss{-\zeta}}
    =
    \frac{\omega^{(\nr-1)\zeta}}{\nr \gauss{\zeta}}
  \end{equation}
\end{proposition}

\begin{proof}
  For \(-1 < \Re\zeta < \nr\), use \(f_\zeta\) and \(I(\zeta)\) from the proof of \cref{thm:exact-quantum-dilogarithm}.
  Since \(f_{\nr-1-\zeta}(w) = -f_\zeta(-w)\),
  \[
    I(\zeta) + I(\nr - 1 - \zeta)
    =
    \int_{\mathbb{R}+\ii\epsilon} f_\zeta(w)\,\mathrm dw
    - \int_{\mathbb{R}-\ii\epsilon} f_\zeta(w)\,\mathrm dw
    =
    -\tu\operatorname{Res}(f_\zeta,0)
  \]
  and we can compute
  \[
    \operatorname{Res}(f_\zeta,0)
    =
    \frac{(2\zeta - \nr + 1)^2}{8\nr}
    - \frac{\nr + 1/\nr}{24}
  \]
  so that
  \[
    I(\zeta) + I(\nr - 1 - \zeta) - 2I(0)
    =
    -\log\nr
    - \frac{\pi\ii}{\nr}\zeta(\zeta - \nr + 1).
  \]
  Exponentiating proves~\eqref{eq:psi-inversion} on the defining strip and we extend by analytic continuation.
\end{proof}

\section{Quantum torus modules and the pentagon relation}
\label{sec:quantum}

\subsection{Quantum torus modules}
\label{sec:quantum-torus-modules}

\begin{definition}
  The \defemph{quantum torus} \(\qtor\) is the algebra over \(\mathbb{C}\) with invertible generators \(\U, \V\) satisfying
  \begin{equation*}
    \U \V = \omega \V \U
  \end{equation*}
\end{definition}

The center of \(\qtor\) is the Laurent polynomial ring in \(\U^{\nr}, \V^{\nr}\), so simple \(\qtor\)-modules are parametrized by central characters \(\chi : \mathbb{C}[\U^{\pm \nr}, \V^{\pm \nr}] \to \mathbb{C}\), which we can identify with \((\mathbb{C} \setminus \set{0})^{2}\) via
\begin{equation*}
  \chi(\U^{\nr}) = a
  \quad
  \chi(\V^{\nr}) = b.
\end{equation*}
The standard way to construct a module with character \(\chi = (a,b)\) is to choose \(\nr\)th roots \(A,B\) of \(a, b\) and use cyclically indexed vectors \(v_{n}, n \in \mathbb{Z}/\nr \mathbb{Z}\) with action
\begin{equation*}
  \U \cdot v_{n} = A \omega^{n} v_{n}
  \quad
  \V \cdot v_{n} = B v_{n+1}.
\end{equation*}
A more equivariant way to do this is to write \(A = \omega^{\alpha}\) and index vectors by choices of \(\alpha\).

\begin{definition}
  For \(a \in \mathbb{C}\) write
  \begin{equation*}
    \Lambda(a) \defeq \set{ \alpha \in \mathbb{C} \given e^{\tu \alpha} = a}
  \end{equation*}
  for the set of logarithms of \(a\).
  For \(a, b \in \mathbb{C} \setminus \set{0}\) consider the \(\nr\)-dimensional space of quasi-periodic functions
  \begin{equation}
    \qpmod{a,b} \defeq
    \set{
      f \colon \Lambda(a) \to \mathbb{C}
      \given
      f(\alpha + \nr) = b^{-1} f(\alpha)
    }
  \end{equation}
  Via the \(\qtor\) action
  \begin{align}
    (\U f)(\alpha) & = \omega^{\alpha} f(\alpha) \\
    (\V f)(\alpha) & = f(\alpha - 1)
  \end{align}
  \(\qpmod{a,b}\) is a simple \(\qtor\)-module with central character \(\chi(\U^{\nr}) = a, \chi(\V^{\nr}) = b\).
\end{definition}
It is natural to view \(\qpmod{} \to \mathcal{X}\) as a rank \(\nr\) vector bundle of \(\qtor\)-modules, where \(\mathcal{X} = ( \mathbb{C} \setminus \set{0} )^{2}\) is the space of characters.
Selecting logarithms \(\alpha, \beta\) is a choice of local coordinates on \(\mathcal{X}\) that allows us to identify the fiber \(\qpmod{e^{\tu \alpha}, e^{\tu \beta}}\) with \(\mathbb{C}^{\nr}\).
One advantage of our description in terms of quasi-periodic functions is that we only need to choose \emph{one} of \(\alpha\) or \(\beta\) to get bases.

If \(\beta \in \Lambda(b)\) then the functions
\begin{equation*}
  \nu(\beta) \colon
  \alpha \mapsto \omega^{-\alpha \beta}
\end{equation*}
are elements of \(\qpmod{a,b}\) and
\begin{equation*}
  \nu(\beta), \dots, \nu(\beta + \nr - 1)
\end{equation*}
form a basis.
We call the \(\nu(\beta)\) the \defemph{frequency basis}, and the vectors satisfy
\begin{equation*}
  \nu(\beta + \nr) = a^{-1} \nu(\beta).
\end{equation*}
They are Fourier dual to the \defemph{delta basis}
\begin{equation*}
  \delta(\alpha)(\alpha')
  =
  \begin{cases}
    b^{(\alpha-\alpha')/\nr} & \text{if } \alpha' - \alpha \in \nr \mathbb{Z} \\
    0                        & \text{otherwise.}
  \end{cases}
\end{equation*}
These vectors satisfy
\begin{equation*}
  \delta(\alpha + \nr) = b \delta(\alpha).
\end{equation*}
The action of \(\qtor\) on these bases is given by
\begin{align}
  \label{eq:qtor-action-frequency} \U \nu(\beta) & = \nu(\beta - 1)                 & \V \nu(\beta)     & = \omega^{\beta} \nu(\beta) \\
  \label{eq:qtor-action-delta} \U \delta(\alpha) & = \omega^{\alpha} \delta(\alpha) & \V \delta(\alpha) & = \delta(\alpha + 1)
\end{align}

\subsection{Maps between modules}
\label{sec:mutations}

We are interested in operators between different modules \(\qpmod{a,b} \to \qpmod{a', b'}\), ie between different fibers of our vector bundle.
In particular we want to understand the \defemph{mutation} operators given by quantum dilogarithms with operator arguments.
Because our dilogarithms have logarithmic arguments we need to take logarithms of the quantum torus generators.
For example the operator
\begin{equation*}
  (\Q f)(\zeta) = \zeta f(\zeta)
\end{equation*}
acts on holomorphic functions of \(\zeta\) and in a suitable sense \(\U = \omega^{\Q}\).
The operator \(\Q\) does not map \(\qpmod{a,b}\) to another space \(\qpmod{a,b}\), but it makes sense to define
\begin{equation*}
  (\psif{\Q} f)(\zeta)
  \defeq
  \psif{\zeta} f(\zeta)
  .
\end{equation*}
If \(f \in \qpmod{a,b}\) then the recurrence relation~\eqref{eq:psi-recurrence} for \(\psif{}\) gives
\begin{equation*}
  (\psif{\Q} f)(\alpha + \nr)
  =
  b^{-1}\psif{\alpha + \nr} f(\alpha)
  =
  b^{-1}(1 - a)^{-1} (\psif{\Q} f)(\alpha)
\end{equation*}
so \(\psif{\Q}\) is an operator
\begin{equation*}
  \psif{\Q}
  \colon
  \mathcal{E}(a,b)
  \to
  \mathcal{E}(a,b(1 - a))
  .
\end{equation*}
It is \emph{not} a map of \(\qtor\)-modules and cannot be one because it changes the central character.
Instead the same recurrence relation gives
\begin{equation*}
  \label{eq:psi-shift-intertwining}
  \V\psif{\Q} = (1 - \U)\psif{\Q}\V.
\end{equation*}
To get the right operators for our pentagon relation we need to add a shift:

\begin{lemma}
  Recall \(\pr = (\nr - 1)/2\) and \(\epsilon = \omega^{\nr \pr} = (-1)^{\nr-1}\).
  The operator
  \begin{equation}
    \psif{\Q + p}
    \colon
    \mathcal{E}(a,b)
    \to
    \mathcal{E}(a,b(1 - \epsilon a))
    \qquad
    (\psif{\Q + p} f)(\alpha)
    \defeq
    \psif{\alpha + p} f(\alpha)
  \end{equation}
  satisfies
  \begin{equation}
    \label{eq:shifted-psi-shift-intertwining}
    \V\psif{\Q + \pr} = (1 - \omega^\pr\U)\psif{\Q + \pr}\V
  \end{equation}
  and commutes with \(\U\).
\end{lemma}

The quantum pentagon relation involves operators like \(\psif{\Q}\), \(\psif{\mathsf{P}}\), and \(\psif{\Q + \mathsf{P}}\) for a suitable logarithm \(\mathsf{P}\) of \(\V\).
Because \(\mathsf{P}\) is a derivative it is difficult to understand its action on any particular \(\qpmod{a,b}\), which consists of functions supported on a discrete domain.
One could consider spaces of generalized functions (eg the Dirac delta function) but this would be an unnecessary digression into functional analysis.
Instead we use Fourier transforms that swap \(\U\) and \(\V\).
Some straightforward computations establish:

\begin{lemma}
  The \defemph{quasi-periodic Fourier transform}
  \begin{equation*}
    \mathcal{F}(f)(\beta) \defeq
    \sum_{k=0}^{\nr-1} \omega^{\beta(\alpha+k)} f(\alpha + k)
  \end{equation*}
  is a linear map
  \begin{equation*}
    \mathcal{F} : \mathcal{E}(a,b) \to \mathcal{E}(b,a^{-1})
  \end{equation*}
  with inverse
  \begin{equation*}
    \mathcal{F}^{-1} : \mathcal{E}(a,b) \to \mathcal{E}(b^{-1},a)
  \end{equation*}
  \begin{equation*}
    \mathcal{F}^{-1}(f)(-\beta)
    \defeq
    \frac{1}{\nr}
    \sum_{k=0}^{\nr-1} \omega^{\beta(\alpha+k)} f(\alpha + k).
  \end{equation*}
  It is not a map of \(\qtor\)-modules, but instead satisfies
  \begin{equation*}
    \mathcal{F} \U \mathcal{F}^{-1} = \V^{-1} , \mathcal{F} \V \mathcal{F}^{-1} = \U.
    \qedhere
  \end{equation*}
\end{lemma}

In summary we have shown:
\begin{proposition}
  The \defemph{mutation} operator \(\K \defeq \mathcal{F}\psif{\Q + p}\) is a linear isomorphism
  \begin{equation*}
    \K \colon
    \mathcal{E}(a,b)
    \to
    \mathcal{E}(b (1 - \epsilon a), a^{-1})
    .
  \end{equation*}
  Writing \(\mu(x) = \K x \K^{-1}\) for conjugation by \(\K\) we have
  \begin{equation*}
    \mu(\U) = \V^{-1}
    \qquad
    \mu(\V) = (1 - \omega^p\V^{-1})^{-1}\U
    .
    \qedhere
  \end{equation*}
\end{proposition}

\subsection{The quantum pentagon relation}

Under repeated applications of \(\K\) the central characters satisfy the same relations~\eqref{eq:five-term-classical-shapes} as in the classical five-term relation.
Write \((a_j,b_j)\) for the central characters after \(j\) applications of \(\K\), starting with \((a_0,b_0) = (a,b)\).
These satisfy
\begin{align*}
  a_{j+1} & = b_j(1 - \epsilon a_j) \\
  b_{j+1} & = a_j^{-1}
\end{align*}
which looks like \cref{eq:five-term-classical-shapes} if we rewrite it as
\begin{equation}
  \label{eq:pentagon-shape-relations}
  1 - a_{j-1} a_{j+1} = \epsilon a_{j}
  \text{ for }
  j \in \mathbb{Z}/5\mathbb{Z}.
\end{equation}
In particular \(a_{5} = a_{0}, b_{5} = b_{0}\) so
\begin{equation*}
  \K^{5} : \qpmod{a,b} \to \qpmod{a,b}
\end{equation*}
as a linear map.
The quantum pentagon relation is the stronger statement that \(\K^{5}\) is a \(\qtor\)-intertwiner independent of \(a\) and \(b\).

\begin{theorem}
  \label{thm:weyl-pentagon-mutation}
  For generic \(a, b\)
  \begin{equation*}
    \K^{5}
    \colon
    \qpmod{a,b} \to \qpmod{a,b}
  \end{equation*}
  is a map of \(\qtor\)-modules and acts by a scalar independent of \(a\) and \(b\):
  \begin{equation*}
    \K^{5}
    = \nr c_{\nr}^{2}\omega^{\pr^{2}}
    = \exp\leftfun(
    \frac{\pi\ii}{6}\left( 4\nr - 9 + \frac{4}{\nr} \right)
    \rightfun).
  \end{equation*}
  Here \(c_{\nr}\) is the constant of \cref{thm:single-qlog-fourier}.
\end{theorem}

\begin{proof}
  We first show that \(\K^{5}\) is a \(\qtor\)-intertwiner; because \(\qpmod{a,b}\) is simple this will show it acts by a scalar.
  To do this we show that \(\mu^{5}\) preserves \(\U\) and \(\V\).
  Set \(X_j = \mu^j(\U)\) so \(X_0 = \U\), \(X_1 = \V^{-1}\), and more generally
  \begin{equation*}
    X_{j+2} = X_j^{-1}(1 - \omega^p X_{j+1}).
  \end{equation*}
  \(\mu\) preserves the quantum torus relation in the sense that
  \begin{equation*}
    X_jX_{j+1} = \omega^{-1}X_{j+1}X_j
  \end{equation*}
  and using this relation to simplify the first one gives
  \begin{align*}
    X_3 & = X_0^{-1} + X_1^{-1} - \omega^{-p}X_0^{-1}X_1^{-1},                               \\
    X_4 & = (1 - \omega^p X_0)X_1^{-1},                        & X_5 & = X_0, & X_6 & = X_1.
  \end{align*}
  In particular \(X_{5} = X_{0}\) and \(X_{6} = X_{1}\), ie
  \begin{equation*}
    \mu^{5}(\U) = \U \text{ and } \mu^{5}(\V^{-1}) = \V^{-1}
  \end{equation*}
  as claimed.
  Since \(\qpmod{a,b}\) is an irreducible \(\qtor\)-module Schur's Lemma tells us \(\K^{5}\) must act by a scalar \(C(a,b)\).

  To show that it is constant we first compute the determinant of \(\K^{5}\).
  Factor \(\K\) into the composition
  \begin{equation*}
    \begin{tikzcd}
      \qpmod{a_{i},b_{i}} \arrow[r, "\psif{\Q + \pr}"] & \qpmod{a_{i},a_{i + 1}} \arrow[r, "\mathcal{F}"] & \qpmod{a_{i + 1}, b_{i+1}}
    \end{tikzcd}
  \end{equation*}
  To compute determinants of the intermediate maps%
  \footnote{Since these are not endomorphisms their determinants depend on a choice of basis; we make compatible choices to obtain the determinant of \(\K^{5}\).}
  choose logarithms \(\alpha_{i}\) of the \(a_{i}\) and write \(\zeta_i = \alpha_i + \pr\).
  By~\eqref{eq:qtor-action-delta} the basis \(\delta(\alpha_{i}), \dots, \delta(\alpha_{i} + \nr - 1)\) of \(\qpmod{a_{i}, b_{i}}\) is a \(\U\)-eigenbasis, so the matrix of \(\psif{\Q + \pr}\) is diagonal and
  \begin{align*}
    \det \psif{\Q + \pr} & = \prod_{k=0}^{\nr-1} \psif{\zeta_i + k}                                                       \\
                         & = \psif{\zeta_i}^{\nr} \prod_{j=1}^{\nr-1}(1 - \omega^{\zeta_i+j})^{-(\nr-j)}                  \\
                         & = \nr^{-\nr/2} \exp\leftfun( -\ldil{\zeta_i} + \frac{\pi\ii}{12}(\nr - 1)(\nr - 2) \rightfun).
  \end{align*}
  Here we used \cref{eq:exact-psi-value} to simplify the product.

  The matrix of \(\mathcal{F}\) is a shifted version of the Vandermonde matrix with entries
  \begin{equation*}
    \left[ \mathcal{F} \right]_{k, l}
    =
    \omega^{(\alpha_{i+1}+k)(\alpha_i+l)},
    \qquad 0 \le k,l < \nr.
  \end{equation*}
  Factoring \(\omega^{\alpha_i(\alpha_{i+1}+k)}\) from row \(k\) and \(\omega^{l\alpha_{i+1}}\) from column \(l\) leaves \(W = (\omega^{kl})_{k,l}\).
  Its Vandermonde determinant is\footnote{Each difference is \(2\ii e^{\pi\ii(k+l)/\nr}\sin(\pi(l - k)/\nr)\).
  The sines are positive so the phases add and \(WW^\dagger = \nr I\) fixes the magnitude.}
  \begin{equation*}
    \det W
    = \prod_{0 \le k<l<\nr}(\omega^l - \omega^k)
    = \nr^{\nr/2}
    e^{\pi\ii(\nr-1)(3\nr-2)/4}.
  \end{equation*}
  We conclude that
  \begin{equation*}
    \det \mathcal{F}
    =
    \nr^{\nr / 2}
    e^{\pi \ii(\nr-1)(3\nr-2) / 4}
    \exp\leftfun(\tu(\zeta_i\zeta_{i+1} - \pr^2)\rightfun)
  \end{equation*}
  and therefore
  \begin{equation*}
    \det \K
    =
    \exp\leftfun(
    \frac{\pi \ii}{6}(\nr - 1)(2\nr - 1)
    + \tu\zeta_i\zeta_{i+1}
    - \ldil{\zeta_i}
    \rightfun)
  \end{equation*}
  so that
  \begin{equation*}
    \det \K^{5}
    =
    \exp\leftfun(
    \frac{5 \pi \ii}{6}(\nr - 1)(2\nr - 1)
    + \sum_{i=0}^{4}\left[ \tu\zeta_i\zeta_{i+1} - \ldil{\zeta_i} \right]
    \rightfun)
  \end{equation*}
  The parameters \(z_i = \epsilon a_i = e^{\tu\zeta_i}\) satisfy~\eqref{eq:five-term-classical-shapes}.
  Applying the classical five-term relation~\eqref{eq:five-term-classical-cyclic} gives
  \begin{equation*}
    C(a,b)^{\nr} = \det \K^{5}
    =
    \exp\leftfun(
    \frac{5 \pi \ii}{6}(\nr - 1)(2\nr - 1)
    - \frac{\pi \ii}{6}
    \rightfun).
  \end{equation*}
  \(C(a,b)\) is a meromorphic function of \((a,b)\) whose \(\nr\)th power \(\det \K^{5} = C(a,b)^{\nr}\) is constant so \(C(a,b)\) is locally constant.
  As in the proof of \cref{thm:five-term-classical} the space of parameters is connected and we conclude \(C(a,b) = C_{\nr}\) is constant.

  To compute the constant \(C_{\nr}\) we take a boundary limit as before.
  For \(0 < t < 1\), choose shapes
  \begin{equation*}
    (z_0,z_1,z_2,z_3,z_4)
    = \leftfun(
    1 - t^2,t,\frac{1}{1 + t},\frac{1}{1 + t},t
    \rightfun).
  \end{equation*}
  These satisfy the classical shape relation~\eqref{eq:five-term-classical-shapes}, namely \(1 - z_{i} = z_{i-1} z_{i+1}\), so the shifted parameters \(a_{i} = \epsilon z_{i}\) satisfy the relation~\eqref{eq:pentagon-shape-relations} with extra signs \(\epsilon\).
  To choose their logarithms set
  \begin{equation*}
    \tau = \frac{\log t}{\tu},
    \qquad
    \mu = -\frac{\log(1 + t)}{\tu},
    \qquad
    \nu = \frac{\log(1 - t^2)}{\tu},
  \end{equation*}
  and choose
  \begin{equation*}
    (\zeta_0,\zeta_1,\zeta_2,\zeta_3,\zeta_4)
    = (\nu + \nr - 1,\tau ,\mu ,\mu ,\tau )
  \end{equation*}
  so that
  \begin{equation*}
    (\alpha_0,\alpha_1,\alpha_2,\alpha_3,\alpha_4)
    = (\nu + \pr,\tau - \pr,\mu - \pr,\mu - \pr,\tau - \pr).
  \end{equation*}
  In the limit \(t \to 0^{+}\) we have \(\tau \to \ii \infty\) and \(\mu, \nu \to 0\).

  Let \(\K_i\) be the matrix of the \(i\)th mutation in the corresponding bases.
  Its entries are
  \begin{equation*}
    (\K_i)_{j,k}
    = \psif{\alpha_i+\pr+k}
    \omega^{(\alpha_i+k)(\alpha_{i+1}+j)}.
  \end{equation*}
  Setting \(k_0 = k_5 = 0\) and \(\alpha_5 = \alpha_0\), we have
  \begin{equation*}
    C_{\nr}
    = (\K_4\K_3\K_2\K_1\K_0)_{0,0}
    = \sum_{k_1,\ldots,k_4=0}^{\nr-1}
    \left[ \prod_{i=0}^{4}\psif{\alpha_i+\pr+k_i} \right]
    \omega^{E(k,t)},
  \end{equation*}
  where the Fourier exponent is
  \begin{align*}
    E(k,t) & = (\tau - \pr + k_1)(\nu + \mu + k_2)       \\
           & \quad + (\tau - \pr + k_4)(\nu + \mu + k_3) \\
           & \quad + (\mu - \pr + k_2)(\mu - \pr + k_3).
  \end{align*}
  For each fixed \(k\) we have \(\ldil{\tau + k} \to \pi\ii/12\) as \(t \to 0^+\) and all logarithmic terms in~\eqref{eq:exact-psi-value} tend to zero, so \(\psif{\tau + k} \to c_{\nr}\).
  Thus the dilogarithm terms converge to
  \begin{equation*}
    \psif{\nr -1} c_{\nr} \psif{k_{2}} \psif{k_{3}} c_{\nr}
    = \frac{c_{\nr}^2}{\nr} \psif{k_{2}} \psif{k_{3}}
  \end{equation*}
  Since \(\tau(\nu + \mu) \to 0\) the only terms in the Fourier factor that survive are those with \(k_{2} = k_{3} = 0\) and we have \(\omega^{E(k,t)} \to \omega^{\pr^2}\).
  The indices \(k_1,k_4\) are free, so the \(\nr^2\) surviving terms give
  \begin{equation*}
    C_{\nr}
    = \nr^2\frac{c_{\nr}^2}{\nr}\omega^{\pr^2}
    = \nr c_{\nr}^2\omega^{\pr^2}.
    \qedhere
  \end{equation*}
\end{proof}

\begin{corollary}
  \label{thm:five-psi-summation}
  For any \(a_j\) satisfying~\eqref{eq:pentagon-shape-relations} and any \(\alpha_{j}\) with \(e^{\tu \alpha_{j}} = a_{j}\),
  \begin{equation}
    \label{eq:five-psi-summation}
    \begin{aligned}
      & \psif{\alpha_0+\pr}\psif{\alpha_3+\pr} \sum_{j=0}^{\nr-1} \psif{\alpha_4+\pr+j} \omega^{(\alpha_4+j)(\alpha_0+\alpha_3)}                                                                                      \\
      & \qquad = \frac{c_{\nr}^{2} \omega^{\pr^{2}}}{\nr^2} \sum_{m,n=0}^{\nr-1} \frac{\omega^{-(\alpha_0)(\alpha_1+m)-(\alpha_1+m)(\alpha_2+n)-(\alpha_2+n)(\alpha_3)}}{\psif{\alpha_1+\pr+m}\psif{\alpha_2+\pr+n}}.
    \end{aligned}
  \end{equation}
\end{corollary}

\begin{proof}
  With respect to the bases \((\delta(\alpha_{i} + k))_{k=0}^{\nr-1}\) the matrices of the mutation operators are
  \begin{align}
    \label{eq:pentagon-kernel} (\K_i)_{j,k}              & = \psif{\alpha_i+\pr+k}\omega^{(\alpha_i+k)(\alpha_{i+1}+j)}, \\
    \label{eq:pentagon-inverse-kernel} (\K_i^{-1})_{k,j} & = \frac{\omega^{-(\alpha_i+k)(\alpha_{i+1}+j)}}
              {\nr\psif{\alpha_i+\pr+k}}
  \end{align}
  and by \cref{thm:weyl-pentagon-mutation} \(\K_4\K_3 = C_{\nr}\K_0^{-1}\K_1^{-1}\K_2^{-1}\).
  Take the \((0,0)\) matrix coefficient and multiply by \(\psif{\alpha_0+\pr}\).
\end{proof}

\subsection{Recovering the summation identities}

\begin{proof}[Proof of~\zcref{thm:single-qlog-fourier}]
  This is a limiting case of \cref{thm:five-psi-summation}.
  For \(0 < A < B < 1\), choose shapes and logarithms
  \begin{equation*}
    (z_0,z_1,z_2,z_3,z_4)
    = \leftfun(
    \frac{A}{B},1 - A,B,
    \frac{1 - B}{1 - A},\frac{B - A}{B(1 - A)}
    \rightfun),
    \qquad
    \alpha_i = \frac{\log z_i}{\tu} - \pr.
  \end{equation*}
  These satisfy~\eqref{eq:five-term-classical-shapes}, so \(a_i = e^{\tu\alpha_i} = \epsilon z_i\) satisfy~\eqref{eq:pentagon-shape-relations}.
  Put
  \begin{equation*}
    \beta = \frac{\log B}{\tu},
    \qquad
    \gamma = - \frac{\log(1 - B)}{\tu}.
  \end{equation*}
  Then for \(B\) fixed and \(A \to 0^+\) we have
  \begin{equation*}
    \alpha_1,\alpha_4 \to -\pr,
    \qquad
    \alpha_2 = \beta - \pr,
    \qquad
    \alpha_3 \to -\gamma - \pr,
    \qquad
    \alpha_0 \to \ii\infty.
  \end{equation*}
  Multiply both sides of~\eqref{eq:five-psi-summation} by \(u = \omega^{\pr \alpha_0}\) to get
  \begin{equation*}
    \begin{aligned}
      & \psif{\alpha_0+p}\psif{\alpha_3+p} \sum_{j=0}^{\nr-1} \psif{\alpha_4+p+j} \omega^{\alpha_{3}(\alpha_4+j)+\alpha_{0}(\alpha_{4}+\pr+j)}                                                                      \\
      & \qquad = \frac{c_{\nr}^{2} \omega^{\pr^{2}}}{\nr^2} \sum_{m,n=0}^{\nr-1} \frac{\omega^{-\alpha_0(\alpha_1+m-\pr)-(\alpha_1+m)(\alpha_2+n)-(\alpha_2+n)(\alpha_3)}}{\psif{\alpha_1+p+m}\psif{\alpha_2+p+n}}.
    \end{aligned}
  \end{equation*}
  In our limit \(u \to 0\) so the only surviving terms are those with \(j = 0\) and \(m = 2 \pr = \nr - 1\).
  Using \(\psif{\alpha_0+\pr} \to c_{\nr}\), \(\psif{0} = 1\), and \(\psif{\nr - 1} = 1/\nr\), we obtain
  \begin{equation*}
    c_{\nr}\psif{-\gamma}\omega^{\pr^2+\pr\gamma}
    = \frac{c_{\nr}^2\omega^{\pr^2}}{\nr}\omega^{-\pr\gamma}
    \sum_{n=0}^{\nr-1}
    \frac{\omega^{\gamma(\beta+n)}}{\psif{\beta+n}},
  \end{equation*}
  or, after rearranging
  \begin{equation}
    \label{eq:pentagon-boundary-inverse-fourier}
    \sum_{n=0}^{\nr-1}
    \frac{\omega^{\gamma(\beta+n)}}{\psif{\beta+n}}
    = \frac{\nr}{c_{\nr}}
    \omega^{(\nr-1)\gamma}\psif{-\gamma}
  \end{equation}
  which is the first identity.

  For the second, apply the first with parameters
  \begin{equation*}
    -\gamma + \nr - 1,
    \qquad
    - \beta - j
  \end{equation*}
  and take a Fourier transform.
  In more detail, this means we multiply by \((c_{\nr}/\nr)\omega^{(\nr-1-\gamma)(\beta+j)}\) and sum over \(0 \le j < \nr\) to obtain
  \begin{align*}
    \sum_{j=0}^{\nr-1}\psif{\beta+j}\omega^{-\gamma(\beta+j)} & = \frac{c_{\nr}}{\nr} \sum_{k=0}^{\nr-1} \frac{\omega^{-k\beta}}{\psif{-\gamma+\nr-1+k}} \sum_{j=0}^{\nr-1}\omega^{-jk} \\
                                                              & = \frac{c_{\nr}}{\psif{-\gamma+\nr-1}}.
  \end{align*}
  The last equality is discrete Fourier orthogonality.
\end{proof}

\begin{proof}[Proof of~\zcref{thm:1psi1-summation}]
  Applying the Fourier transform identities of~\cref{thm:single-qlog-fourier} to parts of~\eqref{eq:five-psi-summation} gives
  \begin{align*}
    \sum_{j=0}^{\nr-1} \psif{\alpha_4+\pr+j}\omega^{(\alpha_4+j)(\alpha_0+\alpha_3)} & = \frac{c_\nr\omega^{-\pr(\alpha_0+\alpha_3)}}
    {\psif{\alpha_0+\alpha_3+\nr-1}}, \\
    \sum_{m=0}^{\nr-1} \frac{\omega^{-(\alpha_1+m)(\alpha_0+\alpha_2+n)}}
    {\psif{\alpha_1+\pr+m}} & = \frac{\nr}{c_\nr} \omega^{-\pr(\alpha_0+\alpha_2+n)}\psif{\alpha_0+\alpha_2+n}.
  \end{align*}
  These are valid by~\eqref{eq:pentagon-shape-relations}:
  \begin{equation*}
    1 - e^{\tu(\alpha_4+\pr)}
    = e^{\tu(\alpha_0+\alpha_3)},
    \qquad
    1 - e^{\tu(\alpha_1+\pr)}
    = e^{\tu(\alpha_0+\alpha_2+n)}.
  \end{equation*}
  Substituting into~\eqref{eq:five-psi-summation} and collecting the factors independent of \(n\) gives
  \begin{equation*}
    \sum_{n=0}^{\nr-1} \frac{\psif{\alpha_0+\alpha_2+n}}{\psif{\alpha_2+\pr+n}} \omega^{-(\alpha_3+\pr)n} \\
    = \nr \omega^{\pr(\alpha_2-\alpha_{3}-\pr)+\alpha_2\alpha_3} \frac{\psif{\alpha_0+\pr}\psif{\alpha_3+\pr}}
      {\psif{\alpha_0+\alpha_3+\nr-1}}.
  \end{equation*}
  This is the desired summation identity up to a change of variables.
  To recover~\cref{eq:1psi1-summation} put
  \begin{equation*}
    a = e^{\tu\alpha},
    \qquad
    b = e^{\tu\beta},
    \qquad
    c = e^{\tu\gamma} = \frac{1 - a}{1 - b},
  \end{equation*}
  and choose
  \begin{equation*}
    \alpha_0 = \alpha - \beta + \pr,
    \qquad
    \alpha_2 = \beta - \pr,
    \qquad
    \alpha_3 = -\gamma - \pr.
  \end{equation*}
  As before write \(\zeta_i = \alpha_i + \pr\) and \(z_i = e^{\tu\zeta_i} = \epsilon a_i\).
  Choose the remaining logarithms \(\alpha_1,\alpha_4\) so that
  \begin{equation*}
    (z_0,z_1,z_2,z_3,z_4)
    = \leftfun(
    \frac{a}{b},
    1 - a,
    b,
    c^{-1},
    1 - \frac{a}{bc}
    \rightfun).
  \end{equation*}
  Now
  \begin{equation*}
    \alpha_0 + \alpha_2 = \alpha,
    \qquad
    \alpha_2 + \pr = \beta,
    \qquad
    - \alpha_3 - \pr = \gamma,
  \end{equation*}
  so the sum on the left is the left-hand side of \cref{eq:1psi1-summation}.
  On the right the exponent simplifies to \((\nr - 1 - \beta)\gamma\) while the three dilogarithm arguments become
  \begin{equation*}
    \alpha_0 + \pr = \alpha - \beta + \nr - 1,
    \qquad
    \alpha_3 + \pr = -\gamma,
    \qquad
    \alpha_0 + \alpha_3 + \nr - 1
    = \alpha - \beta - \gamma + \nr - 1.
  \end{equation*}
  This gives~\eqref{eq:1psi1-summation} and then~\eqref{eq:1psi1-summation-flipped} follows from \cref{eq:psi-inversion,eq:gaussian-key-property}.
\end{proof}

\section{The Fourier transform of the holonomy \texorpdfstring{\(R\)}{R}-matrix}
\label{sec:R-matrix}

In~\cite{arXiv:2509.02354} the author and Reshetikhin found a family \(\holbraid\) of \defemph{holonomy \(R\)-matrices} acting on modules of the de Concini-Kac quantum group \(\mathcal{U}_{\omega^{1/2}}(\mathfrak{sl}_2)\) at a root of unity.
They depend on continuous parameters that can be understood as coordinates for flat connections on the braid complement~\cite{Kashaev2005} or in terms of the geometry of ideal hyperbolic octahedra~\cite{McPhailSnyder2022,McphailSnyder2024} in the \defemph{octahedral decomposition} of the complement.
These matrices are used to define link and tangle invariants (which also depend on the geometric structures) in~\cite{arXiv:2509.02365}.

The operators \(\holbraid\) are defined in terms of their matrix coefficients in the frequency basis \(\nu(\beta)\) of the modules \(\qpmod{a,b}\).
In this section we use \cref{thm:1psi1-summation} to compute the matrix coefficients in the Fourier dual basis \(\delta(\alpha)\).
Topologically this can be understood as passing from the four-term decomposition of the octahedron at the crossing to a five-term decomposition; to do this we apply a series of \(3\)-\(2\) moves on tetrahedra, each corresponding to an application of \cref{thm:1psi1-summation}.
The four-term and five-term decompositions underlie the geometric interpretations of Kashaev's invariant and the colored Jones polynomial developed by \textcite{Yokota2000} and \textcite{Thurston1999}, respectively.
Their relationship is described explicitly by \textcite[Section 2]{ChoMurakami2013}.

In the geometrically degenerate limit where \(g_{1} = g_{2} = (-1)^{\nr-1}\) the braiding \(\holbraid\) recovers the braiding defining the \(\nr\)th colored Jones polynomial at \(q = \omega^{1/2}\).
Taking the matrix coefficients of \(\holbraid\) with respect to the bases \(\nu(\beta)\) recovers Kashaev's \(R\)-matrix~\cite{Kashaev1995}, and we show below that the matrix coefficients for \(\delta(\alpha)\) give the colored Jones \(R\)-matrix.
This gives a new proof of a theorem of \textcite[Proposition 4.2]{Murakami2001}: the Fourier transform of Kashaev's \(R\)-matrix is a specialization of the colored Jones \(R\)-matrix.
Our proof is simpler because we take the Fourier transform \emph{before} taking the geometrically degenerate limit.

\subsection{The space of parameters}

\begin{marginfigure}
  \begin{tikzpicture}[line width=1,baseline=30,scale=2,xscale = 1.5]
    \coordinate (1) at (0,1);
    \coordinate (2) at (0,0);
    \coordinate (3) at (1,0);
    \coordinate (4) at (1,1);
    \draw[->] (2) node[left] {\(2\)} \br (4) node[right] {\(4\)};
    \draw[->] (1) node[left] {\(1\)} \br (3) node[right] {\(3\)};
    \node at (0.5, 1) {\(\lN\)};
    \node at (0, 0.5) {\(\lW\)};
    \node at (0.5, 0) {\(\lS\)};
    \node at (1, 0.5) {\(\lE\)};
  \end{tikzpicture}
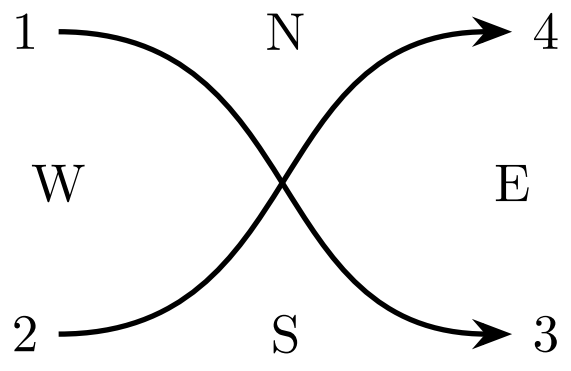
  \caption{Labels for the edges and regions at a crossing.}
  \label{fig:crossing-labels}
\end{marginfigure}

To describe the parameters more succinctly and make the connection to quantum link invariants more explicit we describe them in terms of a braid generator.
Our convention (following~\cite{arXiv:2509.02365}) is to orient it left-to-right and label its edges and regions as in \cref{fig:crossing-labels}.
The holonomy of a flat \(\mathfrak{sl}_{2}\) connection on the complement of a positive braid generator is parametrized by two matrices \(g_{1}, g_{2} \in \slg\) giving the holonomy of the Wirtinger generators and some extra boundary conditions called a \defemph{decoration} and a \defemph{shadow coloring}.%
\note{
  In the language of \textcite{arXiv:1904.10491} these correspond to framings and pinnings, respectively.
}

Concretely these are choices of eigenvectors \(v_{i} \in \mathbb{C}^{2}\) of the \(g_{i}\) and another vector \(u \in \mathbb{C}^{2}\).
These determine labelings of the edges and regions of the diagram as below:
\begin{equation}
  \label{eq:crossing-colored}
  \begin{tikzpicture}[line width=1,baseline=30,scale=2,xscale = 1.5]
    \coordinate (1) at (0,1);
    \coordinate (2) at (0,0);
    \coordinate (3) at (1,0);
    \coordinate (4) at (1,1);
    \draw[->] (2) \br (4);
    \draw[white, line width=10, line cap=butt] (1) \br (3);
    \draw[->] (1) \br (3);
    \node[left] at (1) {\(g_{1}, v_{1}\)};
    \node[left] at (2) {\(g_{2}, v_{2}\)};
    \node[right] at (3) {\(g_{1}, v_{1}\)};
    \node[right] at (4) {\(g_{1}^{-1} g_{2} g_{1}, g_{1}^{-1} v_{2}\)};
    \node at (0.5, 1) {\(u\)};
    \node at (0, 0.5) {\(g_{1} u\)};
    \node at (0.5, 0) {\(g_{2} g_{1} u\)};
    \node at (1, 0.5) {\(g_{1}^{-1} g_{2} g_{1} u\)};
  \end{tikzpicture}
\end{equation}
The rule above is that across an edge with matrix \(g\) the shadow coloring transforms as \(u \mapsto gu\) and that the eigenvectors transform in the compatible way under conjugation.
The \(v_{i}\) determine \defemph{meridian eigenvalues} by
\begin{equation*}
  g_{i} v_{i} = m_{i} v_{i}.
\end{equation*}

Given a decorated, shadow colored representation as above the associated \defemph{octahedral coordinates} are
\begin{align*}
  b_{i} & = \frac{\left\langle \mathbf{e}_2, v_{i} \right\rangle}{\left\langle u_{i^{L}}, v_{i} \right\rangle} & a_{j} & = \left\langle u_{j}, \mathbf{e}_2 \right\rangle
\end{align*}
where \(\left\langle x ,y \right\rangle = \det(x,y)\) is the standard antisymmetric pairing on column vectors, \(\mathbf{e}_2 = (0,1)^{T}\), and \(i^{L}\) is the region to the left of strand \(i\).
(For example \(u_{2^{L}} = g_{1} u\).)
We call the \(b_{i}\) \defemph{segment parameters} and the \(a_{j}\) \defemph{region parameters}.
We call the coloring \(\chi = (g_{1}, g_{2}, v_{1}, v_{2}, u)\) \defemph{admissible} if none of the \(b_{i}\) or \(a_{i}\) is \(0\) or \(\infty\).

These parameters directly determine the dihedral angles of the ideal octahedron assigned to the crossing~\cite{McPhailSnyder2022}.
The four-term decomposition assigns one tetrahedron to each corner of the crossing, and these are assigned shape parameters
\begin{equation}
  \label{eq:crossing-shape-parameters}
  \begin{aligned}
    z_{\lN} & = \frac{b_{4}}{b_{1}}       & z_{\lS} & = \frac{m_{2} b_{2}}{m_{1} b_{3}} \\
    z_{\lW} & = \frac{b_{2}}{m_{1} b_{1}} & z_{\lE} & = \frac{m_{2} b_{4}}{b_{3}}
  \end{aligned}
\end{equation}
We rule out the case \(z_{i} = 0, \infty\) by admissibility above, but it can happen that some \(z_{i} = 1\), in which case we say the crossing is \defemph{pinched}.
In the pinched limit we encounter poles of the quantum dilogarithms \(\psi{}\), but the braiding is still well-defined.
In particular the case \(g_{1} = g_{2} = \epsilon\) recovering Kashaev's \(R\)-matrix is pinched.

\subsection{The holonomy braiding}

We associate to each edge of the crossing~\eqref{eq:crossing-colored} a module \(\qpmod{a,b}\) with central character determined by the octahedral coordinates.
An edge
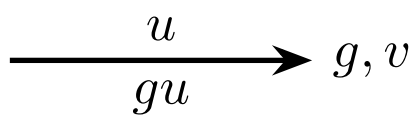
\begin{equation*}
  \begin{tikzpicture}[line width = 1,scale = 2]
    \draw[->] (0,0) to (1, 0) node[right] {\(g, v\)};
    \node[above] at (0.5,0) {\(u\)};
    \node[below] at (0.5,0) {\(g u\)};
  \end{tikzpicture}
\end{equation*}
with parameters
\begin{equation*}
  b =
  \frac{\left\langle \mathbf{e}_{2}, v \right\rangle}{\left\langle u , v \right\rangle}
  ,\quad
  a_{\lN} = \left\langle u, \mathbf{e}_{2} \right\rangle
  ,\quad
  a_{\lS} = \left\langle gu, \mathbf{e}_{2} \right\rangle
\end{equation*}
is assigned the module \(\qpmod{a_{\lS}/a_{\lN}, b}\).
For the diagram~\eqref{eq:crossing-colored} the ratios of region parameters for each edge are
\begin{align*}
  a_{1} & = \frac{a_{\lW}}{a_{\lN}} & a_{4} & = \frac{a_{\lE}}{a_{\lN}} \\
  a_{2} & = \frac{a_{\lS}}{a_{\lW}} & a_{3} & = \frac{a_{\lS}}{a_{\lE}}
\end{align*}
and the holonomy braiding is a linear map
\begin{equation*}
  \holbraid
  =
  \holbraid(\chi)
  \colon
  \qpmod{a_{1}, b_{1}}
  \otimes
  \qpmod{a_{2}, b_{2}}
  \to
  \qpmod{a_{4}, b_{4}}
  \otimes
  \qpmod{a_{3}, b_{3}}
\end{equation*}
or more accurately a family of linear maps depending on the coloring \(\chi\) of the diagram (ie on the geometry of the ideal octahedron).
From a more geometric point of view we have a bundle of vector spaces
\begin{equation*}
  \chi
  \mapsto
  \hom( \qpmod{a_{1}, b_{1}}
  \otimes
  \qpmod{a_{2}, b_{2}}
  ,
  \qpmod{a_{4}, b_{4}}
  \otimes
  \qpmod{a_{3}, b_{3}}
  )
\end{equation*}
over the space of admissible colorings of our diagram and \(\chi \mapsto \holbraid(\chi)\) is a section of this bundle.

Choose logarithms \(\beta_{i}\) of the segment parameters \(b_{i}\) and \(\mu_{k}\) of the meridian eigenvalues \(m_{k}\).
Consider the matrix kernel
\begin{equation}
  \label{eq:holonomy-frequency-kernel}
  \begin{split}
    & \holker{\mu_{1}, \mu_{2}}{\beta_1,\beta_2,\beta_3,\beta_4} \defeq \\
    & \frac{1}{\nr} \frac{\psif{\beta_{4} - \beta_{1}}\psif{\beta_2 - \beta_3 + \mu_2 - \mu_1}}
      {\psif{\beta_4 - \beta_3 + \mu_2}\psift{\beta_2 - \beta_1 - \mu_1}} \omega^{\mu_1(\beta_3-\beta_1)+\mu_2(\beta_2-\beta_4)-\mu_1\mu_2}
  \end{split}
\end{equation}
where we abbreviate
\begin{equation*}
  \psift{\zeta}
  \defeq
  \psif{\zeta + \nr-1}\omega^{(\nr-1)\zeta}.
\end{equation*}
We sometimes write \(\mu\) or \(\beta\) for the corresponding lists of parameters when this does not cause confusion.

\begin{proposition}
  There is a well-defined linear map
  \begin{equation*}
    \mathsf{B} :
    \qpmod{a_{1}, b_{1}}
    \otimes
    \qpmod{a_{2}, b_{2}}
    \to
    \qpmod{a_{4}, b_{4}}
    \otimes
    \qpmod{a_{3}, b_{3}}
  \end{equation*}
  given by
  \begin{equation*}
    \mathsf{B}( \nu(\beta_{1}) \otimes \nu(\beta_{2}) )
    =
    \sum_{n_{3}, n_{4}}
    \holker{\mu}{\beta_{1}, \beta_{2}, \beta_{3} + n_{3}, \beta_{4} + n_{4}}
    \nu(\beta_{4} + n_{4}) \otimes \nu(\beta_{3} + n_{3})
    .
    \qedhere
  \end{equation*}
\end{proposition}

In particular we are going to show that \(\holbraid\) does not depend on the choice of the \(\beta_{i}\).
It \emph{does} depend on the choice of \(\mu_{k}\): these determine the isomorphism classes of the \(\qpmod{a_{i}, b_{i}}\) as quantum group modules~\cite{arXiv:2509.02365}.

\begin{proof}
  Suppose that the crossing is not pinched, so that none of the shape parameters~\eqref{eq:crossing-shape-parameters} is \(1\).
  One can show that the segment and region parameters satisfy the relations
  \begin{equation}
    \label{eq:fourier-module-parameters}
    \begin{aligned}
      a_1 & = m_1\frac{1 - z_{\lW}}{1 - z_{\lN}}, & a_2 & = \frac{1 - z_{\lS}}{m_2(1 - z_{\lW})}, \\
      a_3 & = m_1\frac{1 - z_{\lS}}{1 - z_{\lE}}, & a_4 & = \frac{1 - z_{\lE}}{m_2(1 - z_{\lN})}.
    \end{aligned}
  \end{equation}
  This can be proven by tedious but simple algebraic manipulation.
  More geometrically these identities correspond to the gluing equations of the octahedral decomposition~\cite{McPhailSnyder2022}.
  Here we interpret them as quasi-periodicity conditions.
  For example if we replace \(\beta_{2}\) by \(\beta_{2} + \nr\) so our basis transforms as \(\nu(\beta_{2} + \nr) = a_{2}^{-1} \nu(\beta_{2})\) we also have
  \begin{align*}
    & \holker{\mu}{\beta_{1}, \beta_{2} + \nr, \beta_{3}, \beta_{4}}                                                           \\
    & = \frac{1}{\nr} \frac{\psif{\beta_{4} - \beta_{1} }\psif{\beta_2 + \nr - \beta_3 + \mu_2 - \mu_1}}
        {\psif{\beta_4 - \beta_3 + \mu_2}\psift{\beta_2 + \nr - \beta_1 - \mu_1}} \omega^{\mu_1(\beta_3-\beta_1)+\mu_2(\beta_2+\nr-\beta_4)-\mu_1\mu_2} \\
    & = m_2 \frac{1 - b_{2}/m_{1} b_{1}}{1 - m_{2} b_{2}/m_{1} b_{3}} \holker{\mu}{\beta_{1}, \beta_{2}, \beta_{3}, \beta_{4}} \\
    & = a_2^{-1} \holker{\mu}{\beta_{1}, \beta_{2}, \beta_{3}, \beta_{4}}.
  \end{align*}
  This shows that the operator is independent of the choice of \(\beta_{2}\).
  Similar computations work for the other \(\beta_{i}\).
\end{proof}

\begin{theorem}
  \label{thm:holonomy-fourier}
  The matrix coefficients of the braiding in the delta basis
  \begin{equation*}
    \holbraid ( \delta(\alpha_{1}) \otimes \delta(\alpha_{2}))
    =
    \sum_{n_{3}, n_{4}}
    \holkerhat{\mu}{\alpha_1, \alpha_2, \alpha_3 + n_{3}, \alpha_4 + n_{4}}
    \delta(\alpha_{4} + n_{4}) \otimes \delta(\alpha_{3} + n_{3})
  \end{equation*}
  are \(0\) unless \(\alpha_{1} + \alpha_{2} - \alpha_{3} - \alpha_{4} \in \nr \mathbb{Z}\).
  If \(\alpha_1 + \alpha_2 = \alpha_3 + \alpha_4\) they are
  \begin{equation}
    \label{eq:holonomy-five-psi}
    \begin{aligned}
      \holkerhat{\mu}{\alpha_1, \alpha_2, \alpha_3, \alpha_4} = {} & \gauss{\alpha_2 - \alpha_4} \omega^{\mu_2(\alpha_4-\alpha_2)+\alpha_1\alpha_2}                                                             \\
                                                                   & \times \frac{\psif{\alpha_4 + \mu_2} \psif{\mu_1 - \alpha_3} \psif{\alpha_2 - \alpha_4}}{\psif{\alpha_2 + \mu_2} \psif{\mu_1 - \alpha_1}}.
    \end{aligned}
  \end{equation}
\end{theorem}

\begin{proof}
  Our bases are related by
  \begin{equation*}
    \delta(\alpha) = \frac{1}{\nr}\sum_{n=0}^{\nr-1} \omega^{\alpha(\beta_i+n)}\nu(\beta_i + n)
    \text{ and }
    \nu(\beta) = \sum_{k=0}^{\nr-1} \omega^{-(\beta)(\alpha+k)}\delta(\alpha + k)
  \end{equation*}
  so the matrix coefficients \(\holkerhat{\mu}{}\) are the Fourier transform of the coefficients \(\holker{\mu}{}\):
  \begin{equation}
    \label{eq:holonomy-fourier-four-sum}
    \begin{aligned}
      \holkerhat{\mu}{\alpha} = \frac{1}{\nr^2} \sum_{n_i=0}^{\nr-1} & \holker{\mu}{\beta_1 + n_1, \beta_2 + n_2, \beta_3 + n_3, \beta_4 + n_4}                                \\
                                                                     & \times \omega^{\alpha_1(\beta_1+n_1)+\alpha_2(\beta_2+n_2)-\alpha_3(\beta_3+n_3)-\alpha_4(\beta_4+n_4)}
    \end{aligned}
  \end{equation}
  as suggested by our notation.

  Abbreviate
  \begin{equation}
    \label{eq:fourier-octahedral-logarithms}
    \begin{aligned}
      \zeta_{\lN} & = \beta_4 - \beta_1,                 & \zeta_{\lW} & = \beta_2 - \beta_1 - \mu_1, \\
      \zeta_{\lS} & = \beta_2 - \beta_3 + \mu_2 - \mu_1, & \zeta_{\lE} & = \beta_4 - \beta_3 + \mu_2
    \end{aligned}
  \end{equation}
  and
  \[
    \nu_1 = \alpha_1 - \mu_1,
    \quad \nu_2 = \alpha_2 + \mu_2,
    \quad \nu_3 = \alpha_3 - \mu_1,
    \quad \nu_4 = \alpha_4 + \mu_2
  \]
  and pull out the constant
  \[
    C_0
    =
    \frac{1}{\nr^3}
    \omega^{-\mu_1\mu_2+\beta_1\nu_1+\beta_2\nu_2-\beta_3\nu_3-\beta_4\nu_4-2\pr\zeta_{\lW}}.
  \]
  This shows that~\eqref{eq:holonomy-fourier-four-sum} is equal to
  \begin{equation*}
    C_0\sum_{n_1,n_3=0}^{\nr-1} \omega^{n_1(\nu_1+2\pr)-n_3\nu_3} \psisum{\zeta_{\lN} - n_1}{\zeta_{\lE} - n_3}{-\nu_4} \psisum{\zeta_{\lS} - n_3}{\zeta_{\lW} + 2\pr - n_1}{\nu_2 - 2\pr}
  \end{equation*}
  where we used~\eqref{eq:fourier-module-parameters} to check that
  \[
    e^{-\tu\nu_4}
    = \frac{1 - z_{\lN}}{1 - z_{\lE}},
    \qquad
    e^{\tu(\nu_2-2\pr)}
    = \frac{1 - z_{\lS}}{1 - z_{\lW}}.
  \]
  \Cref{thm:1psi1-summation} and the inversion relation~\eqref{eq:psi-inversion} give
  \begin{align*}
    \psisum{\zeta_{\lN} - n_1}{\zeta_{\lE} - n_3}{-\nu_4}              & = \nr\omega^{(\zeta_{\lE}-n_3)\nu_4}\psif{\nu_4}               \\
                                                                       & \quad \times \frac{\psift{-\zeta + n_3 - n_1}}
    {\psift{-\zeta + n_3 - n_1 + \nu_4}}, \\
    \psisum{\zeta_{\lS} - n_3}{\zeta_{\lW} + 2\pr - n_1}{\nu_2 - 2\pr} & = \nr\omega^{(n_1-\zeta_{\lW})(\nu_2-2\pr)}\psif{2\pr - \nu_2} \\
                                                                       & \quad \times \frac{\psif{\zeta + n_1 - n_3}}
    {\psif{\zeta + n_1 - n_3 - \nu_2 + 2\pr}}
  \end{align*}
  where \(\zeta = \zeta_{\lS} - \zeta_{\lW} = \zeta_{\lE} - \zeta_{\lN} = \beta_1 - \beta_3 + \mu_2\).
  We conclude that~\eqref{eq:holonomy-fourier-four-sum} is equal to
  \begin{equation*}
    \begin{aligned}
      & C_0\nr^2 \omega^{\zeta_{\lE}\nu_4+\zeta_{\lW}(2\pr-\nu_2)} \psif{\nu_4}\psif{2\pr - \nu_2} \\
      & \times\sum_{n_1,n_3=0}^{\nr-1} \omega^{n_1(\nu_1+\nu_2)-n_3(\nu_3+\nu_4)} \frac{\psift{-\zeta + n_3 - n_1}}
      {\psift{-\zeta + n_3 - n_1 + \nu_4}} \frac{\psif{\zeta + n_1 - n_3}}
      {\psif{\zeta + n_1 - n_3 - \nu_2 + 2\pr}}.
    \end{aligned}
  \end{equation*}
  Using
  \begin{equation}
    \label{eq:holonomy-shifted-inversion}
    \psift{-x}
    = \frac{1}{\nr\gauss{x}\psif{x}},
    \qquad
    \frac{\gauss{x - \nu_4}}{\gauss{x}}
    = \gauss{-\nu_4}\omega^{-x\nu_4}
  \end{equation}
  and changing variable to \(n_1 = n_3 + r\) the summand in the double sum becomes
  \[
    \gauss{-\nu_4}\omega^{-\zeta\nu_4}
    \omega^{n_3(\nu_1+\nu_2-\nu_3-\nu_4)}
    \omega^{r(\nu_1+\nu_2-\nu_4)}
    \frac{\psif{\zeta - \nu_4 + r}}{\psif{\zeta - \nu_2 + 2\pr + r}}.
  \]
  The sum over \(n_{3}\) shows that \(\holkerhat{\mu}{\alpha}\) vanishes unless \(\nu_1 + \nu_2 - \nu_3 - \nu_4 = \alpha_{1} + \alpha_{2} - \alpha_{3} - \alpha_{4} \in \nr \mathbb{Z}\) as claimed.
  From now on assume \(\alpha_{1} + \alpha_{2} = \alpha_{3} + \alpha_{4}\) (the general case can be recovered by quasi-periodicity).
  Now the sum over \(n_{3}\) gives a factor \(\nr\) and~\eqref{eq:holonomy-fourier-four-sum} is equal to
  \[
    C_0\nr^3
    \omega^{\zeta_{\lE}\nu_4+\zeta_{\lW}(2\pr-\nu_2)-\zeta\nu_4}
    \psif{\nu_4}\psif{2\pr - \nu_2}\gauss{-\nu_4}
    \psisum{\zeta - \nu_4}{\zeta - \nu_2 + 2\pr}{\nu_3}.
  \]
  Applying~\eqref{eq:1psi1-summation-flipped} results in
  \begin{equation}
    \label{eq:holonomy-fourier-before-inversion}
    \omega^{\Xi}\gauss{-\nu_4}
    \frac{\psif{\nu_4}\psif{2\pr - \nu_2}\psif{2\pr + \nu_3 + \nu_4 - \nu_2}}
    {\psif{2\pr + \nu_3}\psif{2\pr + \nu_4 - \nu_2}},
  \end{equation}
  where all powers of \(\nr\) have cancelled, and
  \begin{align*}
    \Xi & \defeq - \mu_1\mu_2 + \beta_1\nu_1 + \beta_2\nu_2 - \beta_3\nu_3 - \beta_4\nu_4 - 2\pr\zeta_{\lW}   \\
        & \phantom{\defeq} + \zeta_{\lE}\nu_4 + \zeta_{\lW}(2\pr - \nu_2) - \zeta\nu_4 - (\zeta - \nu_4)\nu_3 \\
        & = -\mu_1\mu_2 + \mu_1\nu_2 - \mu_2\nu_3 + \nu_3\nu_4                                                \\
        & = \mu_1(\alpha_2 - \alpha_4) + \alpha_3\alpha_4.
  \end{align*}
  In the second equality the remaining term is \(\beta_1(\nu_1 + \nu_2 - \nu_3 - \nu_4) = 0\).
  Since \(\nu_3 + \nu_4 - \nu_2 = \nu_1\), the four shifted arguments in~\eqref{eq:holonomy-fourier-before-inversion} have the form \(2\pr - x\).
  Applying~\eqref{eq:psi-inversion} in the form
  \[
    \psif{2\pr - x}
    =
    \frac{1}{\nr\gauss{-x}\psif{x}}
  \]
  to all four factors cancels the powers of \(\nr\) and gives
  \[
    \holkerhat{\mu}{\alpha}
    =
    \omega^{\Xi}
    \frac{\gauss{-\nu_4}\gauss{\nu_3}\gauss{\nu_4 - \nu_2}}
    {\gauss{-\nu_2}\gauss{\nu_1}}
    \frac{\psif{\nu_4}\psif{-\nu_3}\psif{\nu_2 - \nu_4}}
    {\psif{\nu_2}\psif{-\nu_1}}.
  \]
  By~\eqref{eq:gaussian-key-property} and \(\nu_{1} + \nu_{2} = \nu_{3} + \nu_{4}\)
  \[
    \frac{\gauss{-\nu_4}\gauss{\nu_3}\gauss{\nu_4 - \nu_2}}
    {\gauss{-\nu_2}\gauss{\nu_1}}
    =
    \gauss{\nu_2 - \nu_4}
    \omega^{(\nu_3-\nu_2)(\nu_2-\nu_4)}.
  \]
  Finally
  \[
    \Xi + (\nu_3 - \nu_2)(\nu_2 - \nu_4)
    =
    \mu_2(\alpha_4 - \alpha_2) + \alpha_1\alpha_2,
  \]
  which proves~\eqref{eq:holonomy-five-psi}.
\end{proof}

\subsection{
  Equivalence of the Kashaev and Jones \texorpdfstring{\(R\)}{R}-matrices
}

Earlier we called the coloring \(\chi\) of our crossing~\eqref{eq:crossing-colored} \defemph{pinched} if any of the shape parameters \(z_{i}\) is \(1\).
Geometrically this is the case where the ideal tetrahedra are degenerate: the shape parameters are the cross-ratios of their ideal vertices, which we can think of as lying in \(\mathbb{C}P^{1}\).
For a pinched crossing the arguments of the dilogarithms in the matrix kernel \(\holker{}{}\) lie in \(\mathbb{Z}\), so they have zeros and poles.
One can check~\cite[Theorem 7.1]{arXiv:2509.02365} that the poles all cancel and the braiding is still well-defined.

Here we consider the case that recovers Kashaev's \(R\)-matrix.
Choose \(g_{1} = g_{2} = \epsilon = (-1)^{\nr-1}\) and \(v_{1} = v_{2}\).
It is not hard to work out that the segment parameters will satisfy
\begin{equation*}
  b_{1} = b_{4} = b
  \text{ and } b_{2} = b_{3} = \epsilon b
\end{equation*}
for some \(b\).
For \(\beta\) with \(e^{\tu \beta} = b\) it is natural to choose
\begin{equation*}
  \beta_{1} = \beta_{4} = \beta \text{ and } \beta_{2} = \beta_{3} = \beta + \pr
\end{equation*}
where \(\pr = (\nr - 1)/2\).
To recover Kashaev's \(R\)-matrix we need to use a rescaled version of the frequency basis
\begin{equation}
  \label{eq:frequency-basis-rescaled}
  e_i(k) \defeq \omega^{\pr k}\nu(\beta_i + k).
\end{equation}
Define \(K\) to be the matrix of the braiding in these bases:
\begin{equation*}
  \holbraid(e_1(k_1) \otimes e_2(k_2))
  =
  \sum_{k_3,k_4=0}^{\nr-1}
  K_{k_1 k_2}^{k_3 k_4}
  e_4(k_4) \otimes e_3(k_3).
\end{equation*}
Equivalently
\begin{equation*}
  K_{k_1 k_2}^{k_3 k_4}
  \defeq
  \omega^{\pr(k_1+k_2-k_3-k_4)}
  \holker{\mu_1,\mu_2}{\beta_1+k_1,\beta_2+k_2,\beta_3+k_3,\beta_4+k_4}
\end{equation*}
where the \(\mu_{i}\) are half-integers with \(e^{\tu \mu_{1}} = e^{\tu \mu_{2}} = \epsilon\).
\begin{proposition}
  When \(\mu_{1} = \mu_{2} = \pr\) we recover Kashaev's \(R\)-matrix~\cite{Kashaev1995}:
  \begin{equation}
    \label{eq:kashaev-cyclic-matrix}
    K_{k_1 k_2}^{k_3 k_4}
    =
    \frac{\omega^{p^2+k_4-k_1}}{\nr}
    \theta_{k_1 k_2}^{k_3 k_4}
    \frac{(\omega)_{[k_4-k_3]}(\omega)_{[k_2-k_1-1]}}
    {(\omega)_{[k_4-k_1]}(\omega)_{[k_2-k_3]}}.
  \end{equation}
  Here
  \[
    (\omega)_n \defeq (\omega;\omega)_n = \prod_{j=1}^{n}(1 - \omega^j),
  \]
  \([r]\) is the residue of \(r\) in \(\{0,\ldots,\nr - 1\}\), and
  \begin{equation*}
    \theta_{k_1 k_2}^{k_3 k_4}
    =
    \begin{cases}
      1 & \text{ if } [k_4 - k_3] + [k_2 - k_1 - 1] - [k_4 - k_1] - [k_2 - k_3] = \nr - 1, \\
      0 & \text{ otherwise.}
    \end{cases}
    \qedhere
  \end{equation*}
\end{proposition}

The matrix~\eqref{eq:kashaev-cyclic-matrix} is related to the matrix \(R_K\) of \textcite[Section 4]{Murakami2001} by
\begin{equation*}
  K_{k_1 k_2}^{k_3 k_4}
  =
  \omega^{\pr^2-1}
  (R_K)_{[-k_3],[-k_4]}^{[-k_2],[-k_1]}.
\end{equation*}
The normalization change is related to framing dependence.

\begin{proof}
  This is a special case of~\cite[Theorem 7.1]{arXiv:2509.02365}, and I suspect some form of this argument was known to \textcite{Kashaev1995} in the early 90s.
  The idea is to keep track of the orders of the zeros of the various dilogarithms, using
  \begin{equation*}
    \psif{\zeta + k}
    =
    \psif{\zeta + [k]}
    (1 - \omega^{\nr \zeta})^{- \gamma(k)}
  \end{equation*}
  where \(\gamma\) is defined by \(k = [k] + \nr\gamma(k)\).
  We conclude that the matrix coefficient has order
  \begin{equation*}
    1 + \gamma(k_4 - k_3) + \gamma(k_2 - k_1 - 1)
    - \gamma(k_4 - k_1) - \gamma(k_2 - k_3)
  \end{equation*}
  which is always nonnegative (so there are no poles) and vanishes exactly when \(\theta_{k_1 k_2}^{k_3 k_4} = 1\).
\end{proof}

With the parameters above we have \(a_{i} = \epsilon\), so it is natural to choose \(\alpha_{i} = \pr\).
To obtain the colored Jones \(R\)-matrix we want to use a rescaled delta basis of \(\qpmod{a_{i},b_{i}}\) given by
\begin{equation*}
  w_{i}(n)
  \defeq
  \nr\omega^{-\beta_{i}(\pr-n)}\delta(\pr - n)
  \text{ for } n = 0, \dots, \nr - 1.
\end{equation*}
(Recall that \(\beta_{1} = \beta_{4} = \beta\) but \(\beta_{2} = \beta_{3} = \beta + \pr\).)
Indexing by \(\pr - n\) corresponds to the usual indexing of weight vectors by their distance from the highest-weight vector \(w_{i}(0)\).
These vectors are Fourier dual to the shifted frequency vectors:
\begin{equation*}
  w_i(n) = \sum_{k=0}^{\nr-1}\omega^{-kn}e_i(k)
  .
\end{equation*}
Write \(J\) for the matrix coefficients of \(\holbraid\) in this basis:
\begin{equation*}
  \holbraid(w_1(n_1) \otimes w_2(n_2))
  =
  \sum_{n_3, n_4=0}^{\nr-1}
  J_{n_1 n_2}^{n_3 n_4}
  w_4(n_4) \otimes w_3(n_3)
  .
\end{equation*}

\begin{proposition}
  \label{thm:colored-jones-rescaled-delta}
  The Fourier dual coefficients recover the colored Jones \(R\)-matrix at \(q = \omega^{1/2}\):
  \begin{equation}
    \label{eq:colored-jones-root-matrix}
    J_{n_1 n_2}^{n_3 n_4}
    =
    \begin{cases}
      \displaystyle \omega^{p^2+n_3(n_2-2p)} \frac{(\omega)_{n_1}(\omega)_{n_4}}
    {(\omega)_{n_2}(\omega)_{n_3}(\omega)_{n_4-n_2}}, & n_1 + n_2 = n_3 + n_4,\quad n_4 \ge n_2, \\[6pt]
      0, & \text{otherwise.}
    \end{cases}
  \end{equation}
\end{proposition}

The matrix~\eqref{eq:colored-jones-root-matrix} is related to the matrix \(R_J\) of \textcite[Section 3]{Murakami2001} by
\begin{equation*}
  J_{n_1 n_2}^{n_3 n_4}
  =
  \omega^{\pr(n_4-n_1)}
  (R_J)_{n_2 n_1}^{n_3 n_4}.
\end{equation*}
The extra change in phase is related to the shift by \(\pr\) in the log-parameters \(\beta_{2} = \beta_{3} = \beta + \pr\) and to the matrix \(D\) in~\cite[Proposition 4.2]{Murakami2001}.

\begin{proof}
  Substitute \(\alpha_i = \pr - n_i\) and \(\mu_1 = \mu_2 = \pr\) into~\eqref{eq:holonomy-five-psi}.
  By \cref{thm:holonomy-fourier} the matrix coefficient~\eqref{eq:holonomy-five-psi} is \(0\) unless
  \begin{equation*}
    n_1 + n_2 - n_3 - n_4 \in \set{-\nr, 0, \nr}
  \end{equation*}
  and it is not hard to show that it vanishes unless \(n_{1} + n_{2} = n_{3} + n_{4}\).
  Thus our matrix coefficient is
  \begin{equation*}
    \gauss{n_{4} - n_{2}}
    \omega^{\pr^2-\pr(n_1+n_4)+n_1n_2}
    \frac{\psif{\nr - 1 - n_4}\psif{n_3}\psif{n_{4} - n_{2}}}
    {\psif{\nr - 1 - n_2}\psif{n_1}}
    .
  \end{equation*}
  Since \(\psif{j} = (\omega)_j^{-1}\) for \(0 \le j < \nr\) we see that \(\psif{n_{4} - n_{2}}\) is \(0\) unless \(n_{4} \ge n_{2}\).
  Furthermore using
  \begin{equation*}
    (\omega)_{\nr-1-j}
    =
    \frac{\nr(-1)^j\omega^{j(j+1)/2}}{(\omega)_j} \text{ and }
    \gauss{k} = (-1)^k\omega^{k(k-1)/2}
  \end{equation*}
  we can work out that
  \begin{equation*}
    \gauss{n_{4} - n_{2}}\frac{(\omega)_{\nr-1-n_2}}{(\omega)_{\nr-1-n_4}}
    =
    \omega^{-(n_2+1)(n_{4}-n_{2})}\frac{(\omega)_{n_4}}{(\omega)_{n_2}}.
  \end{equation*}
  Using \(n_1 = n_3 + t\) and \(2\pr + 1 = \nr\), the unrescaled coefficient therefore becomes
  \begin{equation*}
    \omega^{\pr^2-\pr(n_2+n_3)+n_2n_3}
    \frac{(\omega)_{n_1}(\omega)_{n_4}}{(\omega)_{n_2}(\omega)_{n_3}(\omega)_{n_{4}-n_{2}}}
    .
  \end{equation*}
  Finally the change to the vectors \(w_i(n)\) multiplies this coefficient by
  \begin{equation*}
    \omega^{-\beta_1(\pr-n_1)-\beta_2(\pr-n_2)+\beta_4(\pr-n_4)+\beta_3(\pr-n_3)} = \omega^{\pr(n_2-n_3)}
  \end{equation*}
  which gives~\eqref{eq:colored-jones-root-matrix}.
\end{proof}

\begin{corollary}
  \label{thm:kashaev-colored-jones}
  Write \(F_{jk} = \omega^{-jk}\) for the matrix of the discrete Fourier transform.
  Then
  \begin{equation*}
    J = (F \otimes F)^{-1} K (F \otimes F).\qedhere
  \end{equation*}
\end{corollary}

\begin{proof}
  We are taking the limit where the matrices are both operators
  \begin{equation*}
    \qpmod{\epsilon, b}
    \otimes
    \qpmod{\epsilon, \epsilon b}
    \to
    \qpmod{\epsilon, b}
    \otimes
    \qpmod{\epsilon, \epsilon b}
  \end{equation*}
  In these bases \(F\) is the matrix relating the vectors \(e_i(k)\) and \(w_i(n)\).
\end{proof}

\printbibliography

\end{document}